\documentclass[12pt]{article}
\usepackage{amsfonts}
\usepackage{epsfig}
\usepackage{graphics}
\usepackage{amsmath,amsthm}
\usepackage{amssymb}
\usepackage{inputenc}
\usepackage[T2A]{fontenc}
\begin{document}

\renewcommand{\theequation}{\thesection.\arabic{equation}}
\newtheorem{theorem}{╥хюЁхьр}[section]
\newtheorem{lemma}[theorem]{╦хььр}
\newtheorem{definition}[theorem]{╬яЁхфхыхэшх}
\newtheorem{cor}[theorem]{╤ыхфёЄтшх}
\newtheorem{rem}[theorem]{╟рьхўрэшх}
\newtheorem{example}[theorem]{╧ЁшьхЁ}

\title{╬ эхъюЄюЁ√ї ўрёЄэ√ї ёыєўр ї ЄхюЁхь√ ╨рфюэр-═шъюфшьр фы 
тхъЄюЁэю- ш юяхЁрЄюЁэючэрўэ√ї чрЁ фют}

\author{╤.╤. ┴ющъю, ┬.╩.─єсютющ, └.▀.╒хщЇхЎ}

\maketitle

\section{┬тхфхэшх}

╩ръ шчтхёЄэю, юфэшь шч трцэхщ°шї єЄтхЁцфхэшщ т ЄхюЁшш ьхЁ√ ш
рсёЄЁръЄэюую шэЄхуЁрыр ╦хсхур фы  тх∙хёЄтхээю- ш ъюьяыхъёэючэрўэ√ї
ЇєэъЎшщ  ты хЄё  ЄхюЁхьр ╨рфюэр- ═шъюфшьр (\cite{1} -- \cite{12}).
╧Ёш яюёЄЁюхэшш ЄхюЁшш шэЄхуЁрыр ┴юїэхЁр фы  тхъЄюЁэючэрўэ√ї
ЇєэъЎшщ,
 яЁшэшьр■∙шї чэрўхэш  т
срэрїютюь яЁюёЄЁрэёЄтх, юърч√трхЄё , ўЄю рэрыюушўэр  ЄхюЁхьр
юёЄрхЄё  ёяЁртхфыштющ ыш°№ яЁш эхъюЄюЁ√ї
фюяюыэшЄхы№э√ї яЁхфяюыюцхэш ї. ╬фэръю шчыюцхэшх ¤Єющ ЄхюЁхь√ т
тхъЄюЁэючэрўэюь ёыєўрх т шчтхёЄэ√ї
шёЄюўэшърї ышсю юЄёєЄёЄтєхЄ (\cite{1} -- \cite{6}), ышсю шчырурхЄё 
т ўрёЄэюь ёыєўрх ш ё эхЄюўэюёЄ ьш
(\cite{13}, \cite{15}; ёь. ╩юььхэЄрЁшщ т ъюэЎх ЁрсюЄ√), ышсю фрхЄё 
т ёЄюы№ юс∙хщ ЇюЁьх ёю ёяхЎшЇшўхёъющ
ртЄюЁёъющ ЄхЁьшэюыюушхщ (\cite{7} -- \cite{11}), ўЄю чрЄЁєфэ хЄ
рфхътрЄэюх тюёяЁш Єшх єЄтхЁцфхэш  т
юс∙хяЁшэ Є√ї ЄхЁьшэрї ш юсючэрўхэш ї. ┬ Єю цх тЁхь  фы  трцэ√ї
ёыєўрхт чрЁ фют ёю чэрўхэш ьш т ушы№схЁЄютюь
яЁюёЄЁрэёЄтх шыш т яЁюёЄЁрэёЄтх юуЁрэшўхээ√ї ышэхщэ√ї юяхЁрЄюЁют,
фхщёЄтє■∙шї т ёхярЁрсхы№э√ї ушы№схЁЄют√ї
яЁюёЄЁрэёЄтрї, фюърчрЄхы№ёЄтю ЄхюЁхь√ ╨рфюэр-═шъюфшьр ьюцхЄ с√Є№
яЁютхфхэю т яЁючЁрўэющ ш хёЄхёЄтхээющ ЇюЁьх.
╚чыюцхэш■ ¤Єшї фюърчрЄхы№ёЄт ш яюёт ∙хэр эрёЄю ∙р  ЁрсюЄр. ┼х
яю тыхэшх ёт чрэю Єръцх ё ЁрсюЄющ \cite{14}, уфх
¤Єш єЄтхЁцфхэш  с√ыш ёЇюЁьєышЁютрэ√ ш яЁшьхэхэ√ т ўрёЄэющ ЇюЁьх.
═ю яЁш ¤Єюь т \cite{14} с√ыш фюяє∙хэ√ эхъюЄюЁ√х
эхЄюўэюёЄш, ъюЄюЁ√х ь√ шёяЁрты хь т фрээющ ЁрсюЄх. ┴юыхх яюфЁюсэю
юс ¤Єюь ёърчрэю т ╩юььхэЄрЁшш.

└тЄюЁ√ эшъюшь юсЁрчюь эх яЁхЄхэфє■Є эр юЁшушэры№эюёЄ№
яЁхфырурхь√ї фюърчрЄхы№ёЄт, ъюЄюЁ√х шёяюы№чє■Є шфхш ─рэЇюЁфр ш
╧хЄЄшёр яЁш фюърчрЄхы№ёЄтх тхъЄюЁэюую трЁшрэЄр ЄхюЁхь√
╨рфюэр-═шъюфшьр фы  ёыєўр  срэрїют√ї яЁюёЄЁрэёЄт ё юуЁрэшўхээю
яюыэ√ь срчшёюь ╪рєфхЁр шыш  ты ■∙шїё  ёхярЁрсхы№э√ьш
яЁюёЄЁрэёЄтрьш, ёюяЁ цхээ√ьш ъ фЁєушь срэрїют√ь яЁюёЄЁрэёЄтрь
(яюфЁюсэюх шчыюцхэшх ш шёЄюЁшўхёъшх ётхфхэш  ьюцэю эрщЄш, эряЁшьхЁ, т
\cite{10}). ┬ Єю цх тЁхь , ЁрёёьрЄЁштр  юяхЁрЄюЁэ√щ трЁшрэЄ
ЄхюЁхь√ ╨рфюэр-═шъюфшьр ш єўшЄ√тр  эхёхярЁрсхы№эюёЄ№ яЁюёЄЁрэёЄтр
юуЁрэшўхээ√ї ышэхщэ√ї юяхЁрЄюЁют, фхщёЄтє■∙шї т схёъюэхўэюьхЁэ√ї
ёхярЁрсхы№э√ї ушы№схЁЄют√ї яЁюёЄЁрэёЄтрї, ртЄюЁ√ ёўшЄр■Є
эхюсїюфшь√ь ръЎхэЄшЁютрЄ№ тэшьрэшх ўшЄрЄхы  эр ёє∙хёЄтхээюёЄш
ЄЁхсютрэш  ёхярЁрсхы№эюёЄш ьэюцхёЄтр чэрўхэшщ чрЁ фр ш эр
Ёрчышўш ї, тючэшър■∙шї яЁш ЁрёёьюЄЁхэшш ЁртэюьхЁэющ ш ёшы№эющ
юяхЁрЄюЁэ√ї Єюяюыюушщ (ёь. ╟рьхўрэшх 3.3 ъ ╦хььх 3.2, ╥хюЁхьє 3.5 ш хх
╤ыхфёЄтшх 3.6).

╨рсюЄр шьххЄ ётюхщ Ўхы№■, яЁхцфх тёхую, яЁхфыюцшЄ№ ъЁрЄўрщ°шщ яєЄ№
фы  яюыєўхэш  ¤Єюую трцэюую єЄтхЁцфхэш  схч юсЁр∙хэш 
ъ сюыхх юс∙хщ ЄхюЁшш, р Єръцх юсЁрЄшЄ№ тэшьрэшх ёяхЎшрышёЄют,
ЁрсюЄр■∙шї ё юяхЁрЄюЁэ√ьш ьхЁрьш т ушы№схЁЄют√ї
яЁюёЄЁрэёЄтрї, эр эхъюЄюЁ√х юёюсхээюёЄш ЄхюЁхь√ ╨рфюэр-═шъюфшьр т
¤Єюь ёыєўрх.

\section{╧ЁхфтрЁшЄхы№э√х ётхфхэш }
\setcounter{equation}{0}

┬ ¤Єюь ярЁруЁрЇх ь√ яЁштюфшь эхюсїюфшь√х фы  фры№эхщ°хую ётхфхэш 
юс шчьхЁшь√ї тхъЄюЁэючэрўэ√ї ЇєэъЎш ї ш шэЄхуЁрых ┴юїэхЁр. ┴юыхх
яюфЁюсэє■ шэЇюЁьрЎш■ ьюцэю эрщЄш, эряЁшьхЁ, т \cite{1} --
\cite{4}, \cite{10}, \cite{12}.

╧єёЄ№ $(X, {\mathcal A}, \mu)$ \ - \ яЁюёЄЁрэёЄтю ё ьхЁющ, $E$ \ -
тх∙хёЄтхээюх шыш ъюьяыхъёэюх срэрїютю яЁюёЄЁрэёЄтю,
${\mathcal B} \ = \ {\mathcal B}(E)$ \ - \  $\sigma$-рыухсЁр сюЁхыхтёъшї
ьэюцхёЄт т  $E$. ╘єэъЎш  $f:\, X\to E$ эрч√трхЄё  ёшы№эю
$\mu$-шчьхЁшьющ, хёыш юэр $({\mathcal A}, {\mathcal B})$-шчьхЁшьр
ш $\mu$-ёє∙хёЄтхээю ёхярЁрсхы№эючэрўэр. ╧юёыхфэхх
ючэрўрхЄ, ўЄю ёє∙хёЄтєхЄ Єръюх ьэюцхёЄтю $A\in {\mathcal A}$, ўЄю
$\mu(A)=0$ ш $f(X\setminus A)$  ты хЄё  ёхярЁрсхы№э√ь
яюфьэюцхёЄтюь т $E$. ╘єэъЎш  $f$ эрч√трхЄё  яЁюёЄющ, хёыш юэр $({\mathcal A},
{\mathcal B})$-шчьхЁшьр ш
яЁшэшьрхЄ ъюэхўэюх ўшёыю чэрўхэшщ.

╬ЄьхЄшь, ўЄю шч $({\mathcal A}, {\mathcal B})$-шчьхЁшьюёЄш ЇєэъЎшш
 $f:\, X\to E$ ёыхфєхЄ $\mathcal A$-шчьхЁшьюёЄ№ ЇєэъЎшш $x\to\Vert f(x)\Vert$.
 ─хщёЄтшЄхы№эю, ЇєэъЎш  $h:x\to\Vert x \Vert$  ты хЄё  эхяЁхЁ√тэ√ь юЄюсЁрцхэшхь шч
 $E$ т $\mathbb R$, р, чэрўшЄ, ш ${\mathcal B}$-шчьхЁшьющ. ╧ю¤Єюьє, ёєяхЁяючшЎш 
  $(h\circ f)(x)=\Vert f(x)\Vert$  ты хЄё  $\mathcal A$-шчьхЁшьющ ЇєэъЎшхщ.
\begin{lemma}
╧єёЄ№ $(X, {\mathcal A})$ -- шчьхЁшьюх яЁюёЄЁрэёЄтю ш $E$ --
 срэрїютю яЁюёЄЁрэёЄтю. ╥юуфр фы  $({\mathcal A}, {\mathcal B})$-шчьхЁшьюёЄш ЇєэъЎшш $f:\, X\to E$
 эхюсїюфшью ш фюёЄрЄюўэю, ўЄюс√ фы  ы■сющ эхяЁхЁ√тэющ ЇєэъЎшш $g:E \to \mathbb R$
 ёєяхЁяючшЎш    $g\circ f$ с√ыр ${\mathcal A}$ -- шчьхЁшьющ.
\end{lemma}
{\bf ─юърчрЄхы№ёЄтю.} ═хюсїюфшьюёЄ№ єёыютш  ыхьь√  ты хЄё  юўхтшфэющ. ─ы  фю-ърчрЄхы№ёЄтр
фюёЄрЄюўэюёЄш ЁрёёьюЄЁшь яЁюшчтюы№эюх юЄъЁ√Єюх ьэюцхёЄтю $U\subset E, U \neq E$. ─юёЄрЄюўэю яюърчрЄ№,
ўЄю $f^{-1}(U)\in {\mathcal A}$. ╨рёёьюЄЁшь ЇєэъЎш■ $g_{U} (y)=dist (y,E\backslash U), y \in E$.
╬ўхтшфэю, $g:E \to \mathbb R$-- эхяЁхЁ√тэюх юЄюсЁрцхэшх, яЁш ¤Єюь $U=\{y\in E:g_{U} (y)>0\}$. ╥ръ ъръ
$$
f^{-1}(U)=\{x\in X:(g\circ f)(x)>0\},
$$
Єю $f^{-1}(U)\in {\mathcal A}$.
\begin{theorem} ╧єёЄ№ $(X, {\mathcal A}, \mu)$ -- яЁюёЄЁрэёЄтю ё ьхЁющ, $E$ -- срэрїютю яЁюёЄЁрэёЄтю.
╥юуфр

$(a)$ ьэюцхёЄтю $({\mathcal A}, {\mathcal B})$ шчьхЁшь√ї ЇєэъЎшщ шч $X$ т $E$ чрьъэєЄю юЄэюёшЄхы№эю яюЄюўхўэющ ёїюфшьюёЄш;

$($с$)$ ьэюцхёЄтю ёшы№эю $\mu$-шчьхЁшь√ї ЇєэъЎшщ шч $X$ т $E$ чрьъэєЄю юЄэюёшЄхы№эю яюЄюўхўэющ ёїюфшьюёЄш.
\end{theorem}
{\bf ─юърчрЄхы№ёЄтю.} ╬ўхтшфэю, єЄтхЁцфхэшх $($с$)$ эхяюёЁхфёЄтхээю ёыхфєхЄ шч $(a)$. ─ы  фюърчрЄхы№ёЄтр
єЄтхЁцфхэш  $(a)$ ЁрёёьюЄЁшь яюёыхфютрЄхы№эюёЄ№ $\{f_n(x)\} \ ({\mathcal A}, {\mathcal B})$-шчьхЁшь√ї ЇєэъЎшщ ш
яЁхфяюыюцшь, ўЄю фы  ърцфюую $x\in X$ ёє∙хёЄтєхЄ $\lim\limits_{n\to\infty}f_n(x)$. ╧юърцхь, ўЄю ЇєэъЎш 
$f(x)=\lim\limits_{n\to\infty}f_n(x)$,  ты хЄё  $({\mathcal A}, {\mathcal B})$-шчьхЁшьющ.
─хщёЄтшЄхы№эю, хёыш $g:E \to \mathbb R$-- эхяЁхЁ√тэр  ЇєэъЎш ,Єю
$$
\lim\limits_{n\to\infty}g(f_n(x))=g(f(x)), \ x \in X.
$$
╥хяхЁ№ фюърч√трхьюх єЄтхЁцфхэшх ёыхфєхЄ шч ╦хьь√ 2.1.

┬рцэє■ Ёюы№ яЁш ттхфхэшш шэЄхуЁрыр ┴юїэхЁр шуЁрхЄ
\begin{theorem} ╧єёЄ№ $(X, {\mathcal A}, \mu)$ -- яЁюёЄЁрэёЄтю ё ьхЁющ, $E$ -- срэрїютю яЁюёЄЁрэёЄтю.
╘єэъЎш  $f:\, X\to E$ ёшы№эю $\mu$-
шчьхЁшьр Єюуфр ш Єюы№ъю Єюуфр, ъюуфр юэр
$({\mathcal A}, {\mathcal B})$-шчьхЁшьр ш  ты хЄё  яЁш $\mu$-яюўЄш
тёхї $x\in X$ яЁхфхыюь яюёыхфютрЄхы№эюёЄш яЁюёЄ√ї ЇєэъЎшщ
$f_n:  X\to E, \ n\in {\mathbb N}$. ╧Ёш ¤Єюь яюёыхфютрЄхы№эюёЄ№ $\{f_n\}$
ьюцэю т√сЁрЄ№ Єръ, ўЄюс√ яЁш ы■сюь $n\in {\mathbb N}$ яЁш $\mu$-яюўЄш
тёхї $x\in X$ т√яюыэ ыюё№ эхЁртхэёЄтю $\Vert f_n (x)\Vert
\le \Vert f(x)\Vert$.
\end{theorem}

{\bf ─юърчрЄхы№ёЄтю.} ╬ўхтшфэю т фюърчрЄхы№ёЄтх эєцфрхЄё  Єюы№ъю эхюсїюфшьюёЄ№ єёыютшщ ЄхюЁхь√.
═х эрЁє°р  юс∙эюёЄш, ьюцэю яЁхфяюыюцшЄ№, ўЄю $f(X)$  ты хЄё  ёхярЁрсхы№э√ь ш ёюфхЁцшЄ эхэєыхтющ ¤ыхьхэЄ
шч $E$. ╬сючэрўшь ўхЁхч $K$ эх сюыхх ўхь ёўхЄэюх яюфьэюцхёЄтю т $f(X)$ Єръюх, ўЄю $f(X) \subset \overline{K}$.
╧єёЄ№ ${\mathbb Q}$ \ - \ ьэюцхёЄтю ЁрЎшюэры№э√ї ўшёхы ш
$$
K_{\mathbb Q} : = \{ \ qk \ : \ k\in K , \ q\in {\mathbb Q} \ \}.
$$
╥ръ ъръ $K_{\mathbb Q}$ \ - \ ёўхЄэюх ьэюцхёЄтю, яЁхфёЄртшь хую т тшфх $K_{\mathbb Q} \ = \ \{ y_n \ : \ n\in {\mathbb N}\}$, яЁш ¤Єюь ьюцэю ёўшЄрЄ№, ўЄю $y_1 = 0$.
╧юърцхь, ўЄю фы  ърцфюую $y$ шч $f(X)$ ш яЁюшчтюы№эюую $\varepsilon > 0$ ёє∙хёЄтєхЄ ¤ыхьхэЄ $y_m \in K_{\mathbb Q}$ Єръющ, ўЄю
ёяЁртхфышт√ ёююЄэю°хэш 
\begin{equation}
 \Vert y_m\Vert \le \Vert y\Vert ,  \ \ \ \ \ \  \Vert { y_m - y } \Vert < \varepsilon.  \label{21}
\end{equation}
═х эрЁє°р  юс∙эюёЄш, ьюцэю ёўшЄрЄ№, ўЄю
$ y\neq 0 $ \ ш \ $ 0 < \varepsilon  < \Vert y\Vert $. ┬√схЁхь $\varepsilon_1$ Єръ, ўЄюс√ т√яюыэ ышё№ єёыютш 
\ $ 0 < \varepsilon_1 < 2^{-1} \varepsilon $ \ ш \ ${\Vert y\Vert}^{-1}{\varepsilon_1} \in {\mathbb Q}.$
╧юыюцшь $ \alpha : = 1 - {\Vert y\Vert}^{-1}{\varepsilon_1}.$ ╥юуфр
$$
\alpha\in {\mathbb Q}, \ \  0 < \alpha < 1, \ \ \Vert y\Vert - \varepsilon_1 = \alpha \Vert y\Vert, \ \ \Vert {y - \alpha y } \Vert = \varepsilon_1.
$$
╚ч юяЁхфхыхэш  ьэюцхёЄтр $K$ ёыхфєхЄ, ўЄю ёє∙хёЄтєхЄ ¤ыхьхэЄ $k \in K$, Єръющ, ўЄю
$$
\Vert {y - k} \Vert < \varepsilon_1.
$$
┼ёыш  $\Vert k\Vert \le \Vert y\Vert$, Єю т ърўхёЄтх ¤ыхьхэЄр $y_m$ тюч№ьхь $k$, т яЁюЄштэюь ёыєўрх яюыюцшь $y_m : = \alpha k.$
─хщёЄтшЄхы№эю, Єюуфр шьххь
$$
\Vert y_m \Vert = \Vert \alpha k\Vert = \Vert ( \alpha y - \alpha k ) - \alpha y \Vert \leq \alpha \Vert {y - k} \Vert + \alpha \Vert y\Vert
$$
$$
< \alpha \varepsilon_1 +  \alpha \Vert y\Vert = \alpha\varepsilon_1 + \Vert y\Vert  - \varepsilon_1 = \Vert y\Vert - ( 1 - \alpha)\varepsilon_1 < \Vert y\Vert,                                             $$
яЁш ¤Єюь
$$
\Vert {y_m - y} \Vert = \Vert {\alpha k - y} \Vert = \Vert (\alpha k - \alpha y) + (\alpha y - y) \Vert
$$
$$
\leq \alpha \Vert {k - y} \Vert + \Vert {y - \alpha y} \Vert < \alpha \varepsilon_1 +\varepsilon_1 < 2\varepsilon_1 < \varepsilon.
$$
╥ръшь юсЁрчюь, єёыютш  (2.1) т√яюыэ ■Єё .

%─ы  ърцфюую $x\in X$ ш ърцфюую $n\in {\mathbb N}$ юяЁхфхышь яюфьэюцхёЄтю
%$$
%A_n (x) = \{ y_j : \ y_j \in K_{\mathbb Q}, \ j \leq n, \  \Vert y_j\Vert \le \Vert f(x)\Vert \  \}.
%$$
%╥ръ ъръ $y_1 = 0$, Єю $ A_n (x) \neq \emptyset$ фы  ърцфюую $x\in X$ ш ърцфюую $\mathbb Q$.

╥хяхЁ№ яхЁхщфхь ъ яюёЄЁюхэш■ яюёыхфютрЄхы№эюёЄш $\{f_n\}_{n=1}^{\infty}$. ╘єэъЎшш $f_n$ сєфхь шёърЄ№ т тшфх
\begin{equation}
 f_{n}(x) : = \sum_{j=1}^n y_j \chi_{B_{n,j}}(x), \ \ x \in X, \ \label{22}
\end{equation}
уфх $\{B_{n,j}\}_{j=1}^{n} \subset {\mathcal A}, \ \ B_{n,i} \cap B_{n,j} = \emptyset \ ( i \neq j ), \ \ \bigcup_{j=1}^{n} \ B_{n,j} = X.$
┬√сЁрт $B_{1,1} : = X,$ Єю хёЄ№ юяЁхфхышт $f_{1}(x) : = 0, x \in X,$ фры№эхщ°шх яюёЄЁюхэш  сєфхь яЁютюфшЄ№ яю шэфєъЎшш. ┼ёыш ЇєэъЎш  $f_{n}(x)$
тшфр (2.2) єцх яюёЄЁюхэр, Єю яюыюцшь
$$
B_{n+1,n+1} : = \{x \in X : \Vert {y_{n+1} - f(x)} \Vert < \Vert {f_n(x) - f(x)} \Vert  \ \}  \ \bigcap \ \ \{ \ x\in X : \Vert y_{n+1} \Vert \leq \Vert f(x) \Vert \  \},
$$
$$
B_{n+1,j} : = B_{n,j} \backslash B_{n+1,n+1} \ ( j = 1,2,...,n ) , \ \ f_{n+1}(x) : = \sum_{j=1}^{n+1} y_j \chi_{B_{n+1,j}}(x), \ \ x \in X.
$$
▀ёэю, ўЄю $\{B_{n+1,j}\}_{j=1}^{n+1} \subset {\mathcal A}, \ \ B_{n+1,i} \cap B_{n+1,j} = \emptyset \ ( i \neq j ), \ \ \bigcup_{j=1}^{n+1} \ B_{n+1,j} = X.$
╤ыхфютрЄхы№эю, ърцфр  шч ЇєэъЎшщ $f_n$  ты хЄё  яЁюёЄющ ш яю яюёЄЁюхэш■ єфютыхЄтюЁ хЄ эхЁртхэёЄтє $\Vert f_{n}(x) \Vert \leq \Vert f(x) \Vert$ яЁш тёхї $x \in X.$
╚ч яюёЄЁюхэш  Єръцх тшфэю, ўЄю яЁш тёхї $n\in {\mathbb N}$ ш тёхї $x \in X.$ т√яюыэ хЄё  эхЁртхэёЄтю
$$
 \Vert {f_{n+1}(x) - f(x)} \Vert \leq \Vert {f_n(x) - f(x)} \Vert,
$$
р шч юЎхэюъ (2.1) ёыхфєхЄ, ўЄю яЁш тёхї $x \in X$ шьххЄ ьхёЄю ЁртхэёЄтю
$$
 \inf_{n\in {\mathbb N}} \Vert {f_n(x) - f(x)} \Vert = 0.
$$
╥ръшь юсЁрчюь, $ f(x) = \lim\limits_{n\to\infty} f_n(x) $ яЁш тёхї $\ x \in X,$ ш ЄхюЁхьр фюърчрэр.

╘єэъЎш  $f:\, X\to E$ эрч√трхЄё  $\mu$-шэЄхуЁшЁєхьющ яю ┴юїэхЁє,
хёыш юэр ёшы№эю $\mu$-шчьхЁшьр ш ЇєэъЎш  $\Vert f(x)\Vert$
 ты хЄё  $\mu$-шэЄхуЁшЁєхьющ яю ╦хсхує. ╤хьхщёЄтю $\mu$-шэЄхуЁшЁєхь√ї яю ┴юїэхЁє
 $E$-чэрўэ√ї ЇєэъЎшщ юсючэрўр■Є ${\mathcal L}^1 (X, \mathcal A, \mu, E)$.
 ╧Ёш $E={\mathbb R}$ ($E= {\mathbb C}$) ¤Єю ёхьхщёЄтю  ты хЄё 
 ёхьхщёЄтюь $\mu$-шэЄхуЁшЁєхь√ї яю ╦хсхує ЇєэъЎшщ, юсючэрўрхь√ь
 юс√ўэю
${\mathcal L}^1 (X, \mathcal A, \mu)$ (${\mathcal L}^1 (X, \mathcal A, \mu, {\mathbb C}))$.
╚эЄхуЁры юЄ Єръшї ЇєэъЎшщ юяЁхфхы хЄё  ёыхфє■∙шь юсЁрчюь. ╧Ёхфяюыюцшь тэрўрых,
ўЄю ЇєэъЎш  $f(x)$  ты хЄё  яЁюёЄющ ш $\mu$-шэЄхуЁшЁєхьющ. ╧єёЄ№ $\{a_j\}_{j=1}^n$ - ьэюцхёЄтю эхэєыхт√ї чэрўхэшщ
ЇєэъЎшш $f(x)$, яЁш ¤Єюь $a_j \neq a_k$, хёыш $j \neq k.$ ╧єёЄ№ фрыхх $ A_j$ -  ьэюцхёЄтю, эр ъюЄюЁюь $f(x)$
яЁшэшьрхЄ чэрўхэшх $a_j, \ j = 1,2,...,n.$ ╥юуфр шч шчьхЁшьюёЄш $f(x)$ ёыхфєхЄ шчьхЁшьюёЄ№ ърцфюую $ A_j$, р
шч $\mu$-шэЄхуЁшЁєхьюёЄш $f(x)$ ёыхфєхЄ, ўЄю ърцфюх $ A_j$ шьххЄ ъюэхўэє■ ьхЁє. ╥ръшь юсЁрчюь, т√Ёрцхэшх $\sum_{j=1}^n a_j \mu(A_j)$
шьххЄ ёь√ёы, ш т ¤Єюь ёыєўрх шэЄхуЁры ┴юїэхЁр юяЁхфхы ■Є ЁртхэёЄтюь
$$
\int\limits_X f(x) d\mu(x) \left ( =\int\limits_X fd\mu\right ) :=
\sum_{j=1}^n a_j \mu(A_j).
$$
╦хуъю тшфхЄ№, ўЄю Єюуфр

$(a) \ \ \Vert \ \int\limits_X f(x) d\mu(x) \ \Vert  \ \le \
\sum_{j=1}^n \Vert a_j \Vert \mu(A_j) \ = \ \int\limits_X \Vert f(x) \Vert d\mu(x);$

(с)\ \ эр ьэюцхёЄтх яЁюёЄ√ї  $\mu$-шчьхЁшь√ї ЇєэъЎшщ шэЄхуЁры ┴юїэхЁр  ты хЄё 

\ \ \ \ \ ышэхщэ√ь юяхЁрЄюЁюь.

─ы  яЁюшчтюы№эющ ЇєэъЎшш $f\in {\mathcal L}^1 (X, \mathcal A, \mu, E)$
шэЄхуЁры ┴юїэхЁр юяЁхфхы хЄё  яюёЁхфёЄтюь яЁхфхы№эюую яхЁхїюфр
\begin{equation}
\int\limits_X fd\mu :=\lim\limits_{n\to\infty} \int\limits_X f_n d\mu,  \label{23}
\end{equation}
уфх $\{f_n \}$ -- яЁюшчтюы№эр  яюёыхфютрЄхы№эюёЄ№ яЁюёЄ√ї ЇєэъЎшщ
шч ${\mathcal L}^1 (X, \mathcal A, \mu, E)$, єфютыхЄтюЁ ■∙р  єёыютш ь
╥хюЁхь√ 2.3. ╚ч юЎхэъш
$$
 \Vert \int\limits_X f_n(x) d\mu(x) -  \int\limits_X f_m(x) d\mu(x) \Vert  \ \le \
\int \limits_X \Vert {f_n(x) -  f_m(x)} \Vert d\mu(x)
$$
ш ёююЄэю°хэшщ
$$
\Vert {f_n(x) -  f_m(x)} \Vert \leq \Vert f_n(x) \Vert + \Vert f_m(x) \Vert  \leq  2 \Vert f(x) \Vert,
$$
$$
\lim \limits_{n,m\to\infty}\Vert {f_n(x) -  f_m(x)} \Vert \ =  \ \Vert {f(x) -  f(x))} \Vert = 0,
$$
ъюЄюЁ√х т√яюыэ ■Єё  яЁш  $\mu$-яюўЄш тёхї $x \in X$, шёяюы№чє  ЄхюЁхьє ╦хсхур ю яЁхфхы№эюь яхЁхїюфх ё
ёєььшЁєхьющ ьрцюЁрэЄющ, яюыєўрхь,
ўЄю яЁхфхы т яЁртющ ўрёЄш (2.3) ёє∙хёЄтєхЄ. └эрыюушўэю фюърч√трхЄё , ўЄю ¤ЄюЄ яЁхфхы эх чртшёшЄ юЄ
т√сюЁр яюёыхфютрЄхы№эюёЄш $\{f_n \}$.

╤шы№эю $\mu$-шчьхЁшьр  ЇєэъЎш  $ f : X \rightarrow E $ эрч√трхЄё  $\mu$-шэЄхуЁшЁєхьющ т ёь√ёых ┴юїэхЁр яю
ьэюцхёЄтє $A\in {\mathcal A}$, хёыш $f\cdot \chi_A \in
{\mathcal L}^1 (X, \mathcal A, \mu, E)$, уфх $\chi_A$ -- їрЁръЄхЁшёЄшўхёър 
ЇєэъЎш  ьэюцхёЄтр $A$ ш шэЄхуЁры ┴юїэхЁр ЇєэъЎшш $f$ яю $A$
юяЁхфхы хЄё  ЇюЁьєыющ
$$
\int\limits_A fd\mu :=\int\limits_X f\cdot \chi_A d\mu.
$$

═хЄЁєфэю тшфхЄ№, ўЄю ёхьхщёЄтю ${\mathcal L}^1 (X, \mathcal A, \mu, E)$  ты хЄё  ышэхщэ√ь
яЁюёЄЁрэёЄтюь, р шэЄхуЁры ┴юїэхЁр юсырфрхЄ эр ¤Єюь ёхьхщёЄтх ётющёЄтюь ышэхщэюёЄш, яЁш ¤Єюь
фы  ЇєэъЎшщ шч ${\mathcal L}^1 (X, \mathcal A, \mu, E)$ ёяЁртхфыштр юЎхэър
\begin{equation}
 \Vert \int\limits_X fd\mu  \Vert \le \int\limits_X \Vert f(x)\Vert d\mu(x). \  \label{24}
\end{equation}

%$$
%\Vert f\Vert_1 = \int\limits_X \Vert f(x)\Vert d\mu(x).
%$$
╥ръшь юсЁрчюь, ёхьхщёЄтю ${\mathcal L}^1 (X, \mathcal A, \mu, E)$  ты хЄё  ышэхщэ√ь
яЁюёЄЁрэёЄтюь ё яюыєэюЁьющ
$\Vert f\Vert_1 = \int\limits_X \Vert f(x)\Vert d\mu(x).$ ╧Ёютюф  юс√ўэ√х
т ЄхюЁшш ьхЁ√ Ёрёёєцфхэш , яюыєўрхь, ўЄю ёЄрэфрЁЄэр  яЁюЎхфєЁр юЄюцфхёЄтыхэш  $\mu$-яюўЄш тё■фє ёютярфр■∙шї
ЇєэъЎшщ яЁштюфшЄ ъ срэрїютюьє яЁюёЄЁрэёЄтє, юсючэрўрхьюьє
$L^1 (X, {\mathcal A}, \mu, E)$. ╧Ёш ¤Єюь шэЄхуЁры ┴юїэхЁр
ёЄрэютшЄё  ышэхщэ√ь юяхЁрЄюЁюь, фхщёЄтє■∙шь шч
$L^1 (X, {\mathcal A}, \mu, E)$ т $E$ ш  ты ■∙шьё  ёцрЄшхь (ёь.(2.4)). ╧Ёютюф  Ёрёёєцфхэш , рэрыюушўэ√х Єхь,
ўЄю с√ыш яЁютхфхэ√ т√°х яЁш юсюёэютрэшш ъюЁЁхъЄэюёЄш юяЁхфхыхэш  шэЄхуЁрыр ┴юїэхЁр яюыєўрхь,ўЄю
%$$
%\left \Vert \int\limits_X fd\mu\right \Vert \le \Vert f\Vert_1
%\left ( =\int\limits_X \Vert f(x)\Vert d\mu(x)\right ).
%$$
фы  $E$-чэрўэ√ї ЇєэъЎшщ ёяЁртхфышт рэрыюу ЄхюЁхь√ ╦хсхур ю
яЁхфхы№эюь яхЁхїюфх ё ёєььшЁєхьющ ьрцюЁрэЄющ.

\begin{theorem}
 ╧єёЄ№ $g\in {\mathcal L}^1 (X, A, \mu)$ -- $\mu$-яюўЄш
тё■фє эхюЄЁшЎрЄхы№эр  ЇєэъЎш  ш яєёЄ№ $\{f_n\}$ --
яюёыхфютрЄхы№эюёЄ№ ёшы№эю $\mu$-шчьхЁшь√ї $E$-чэрўэ√ї ЇєэъЎшщ эр
$X$, фы  ъюЄюЁ√ї яЁш $\mu$-яюўЄш тёхї $x\in X$ т√яюыэ ■Єё 
ёююЄэю°хэш :
$$
\Vert f_n (x)\Vert \le g(x), \quad n\in {\mathbb N},
$$
$$
\lim\limits_{n\to\infty} f_n(x) = f(x).
$$
уфх $f$ -- ёшы№эю $\mu$-шчьхЁшьр  ЇєэъЎш .

╥юуфр $\{f,f_1, f_2, \ldots\} \subset {\mathcal L}^1 (X,
{\mathcal A}, \mu, E)$ ш
$$
\int\limits_X f(x) d\mu(x) = \lim\limits_{n\to\infty} \int\limits_X
f_n (x) d\mu(x).
$$
\end{theorem}

╧єёЄ№ $f\in {\mathcal L}^1 (X, {\mathcal A}, \mu, E)$
ш $\varphi \in E^*$ ($E^*$ - яЁюёЄЁрэёЄтю, ёюяЁ цхээюх ъ $E$). ╥юуфр
т ёююЄтхЄёЄтшш ё (2.3) яюыєўрхь
$$
\varphi \left ( \int\limits_X f(x)d\mu(x)\right ) =  \varphi \ (\lim\limits_{n\to\infty}\int\limits_X f_n(x) d\mu(x) )
$$
$$
= \lim\limits_{n\to\infty}\varphi \left ( \int\limits_X f_n (x)d\mu(x)\right ) = \lim\limits_{n\to\infty} \ \int\limits_X (\varphi\circ f_n)(x) d\mu(x).
$$
╥ръ ъръ $\mu$-яюўЄш тё■фє
$$
 \lim\limits_{n\to\infty} \varphi \ ( f_n(x) ) =  \varphi \ ( f(x) ),
$$
$$
 | \  \varphi \ ( f_n(x) ) | \leq \Vert \ \varphi \ \Vert \ \Vert \ f_n(x) \ \Vert \leq \Vert \ \varphi \ \Vert \ \Vert f(x) \ \Vert,
$$
Єю
$$
\varphi \left ( \int\limits_X f(x)d\mu(x)\right ) = \lim\limits_{n\to\infty} \ \int\limits_X (\varphi\circ f_n)(x) d\mu(x) = \int\limits_X (\varphi\circ f)(x) d\mu(x).
$$

╬cюсюх чэрўхэшх шьххЄ ёыєўрщ, ъюуфр $E={\c B}
(H,G)$ - срэрїютю яЁюёЄЁрэёЄтю юуЁрэшўхээ√ї ышэхщэ√ї юяхЁрЄюЁют,
фхщёЄтє■∙шї шч юфэюую ушы№схЁЄютр яЁюёЄЁрэёЄтр $H$ т фЁєуюх ушы№схЁЄютю
яЁюёЄЁрэёЄтю $G$ (хёыш $H=G$, Єю тьхёЄю ${\c B} (H,H)$
єяюЄЁхсы ■Є
ёшьтюы ${\c B} (H)$. ╬яЁхфхыхээ√щ т√°х шэЄхуЁры эрч√тр■Є т
¤Єюь ёыєўрх ЁртэюьхЁэ√ь шэЄхуЁрыюь ┴юїэхЁр, шьх  т тшфє
ЁртэюьхЁэє■ Єюяюыюуш■ т ${\c B} (H,G)$, яюЁюцфрхьє■
юяхЁрЄюЁэющ эюЁьющ, ш юсючэрўр■Є $(u)\int\limits_X fd\mu$.
╘єэъЎшш шч ${\mathcal L}^1 (X, \mathcal A, \mu, {\c B} (H,G))$
эрч√тр■Є ЁртэюьхЁэю шэЄхуЁшЁєхь√ьш яю ┴юїэхЁє. ┬ьхёЄх ё Єхь,
шёяюы№чє  ёшы№эє■ Єюяюыюуш■ т ${\c B} (H,G)$, яюЁюцфрхьє■
яюЄюўхўэющ ёїюфшьюёЄ№■ юяхЁрЄюЁют, ьюцэю ттхёЄш сюыхх юс∙хх
яюэ Єшх шэЄхуЁрыр.

╘єэъЎш■ $f:\, X\to {\c B} (H,G)$ эрч√тр■Є ёшы№эю шэЄхуЁшЁєхьющ
 яю ┴юїэхЁє, хёыш яЁш тёхї $h\in H$ т√яюыэ хЄё  $f(x) h\in
{\mathcal L}^1 (X, \mathcal A, \mu, G)$. ╤хьхщёЄтю Єръшї ЇєэъЎшщ юсЁрчєхЄ
ышэхщэюх ьэюуююсЁрчшх, ёюфхЁцр∙хх т ёхсх $\mathcal L^1 (X, \mathcal A,
\mu, {\c B} (H,G))$.

\begin{theorem}
 ╧єёЄ№ $f:\, X\to {\c B}(H,G)$ -- ёшы№эю
шэЄхуЁшЁєхьр  яю ┴юїэхЁє ЇєэъЎш . ╥юуфр ЇюЁьєыр
\begin{equation}
Th:=\int\limits_X (f(x)h) d\mu(x),\quad h\in H,   \  \label{25}
\end{equation}
юяЁхфхы хЄ юуЁрэшўхээ√щ ышэхщэ√щ юяхЁрЄюЁ $T\in \c B(H,G)$.
\end{theorem}
─юърчрЄхы№ёЄтю. ▀ёэю, ўЄю ЁртхэёЄтюь (2.5) юяЁхфхы хЄё  ышэхщэ√щ юяхЁрЄюЁ $T$, юяЁхфхыхээ√щ
эр тёхь яЁюёЄЁрэёЄтх $H$ ш фхщёЄтє■∙шщ т яЁюёЄЁрэёЄтю $G$.
%╧ю¤Єюьє, фы  фюърчрЄхы№ёЄтр хую
%юуЁрэшўхээюёЄш фюёЄрЄюўэю фюърчрЄ№ хую чрьъэєЄюёЄ№. ─ы  ¤Єюую
╨рёёьюЄЁшь тёяюьюурЄхы№эюх
юЄюсЁрцхэшх $W$ шч яЁюёЄЁрэёЄтр $H$ т $ L^1 (X, \mathcal A,
\mu, G)$, яюыюцшт
$$
Wh \ := \ f(x)h, \ \ h \in H.
$$
╬ўхтшфэю, юЄюсЁрцхэшх $ W $  ты хЄё  ышэхщэ√ь. ─ы  фюърчрЄхы№ёЄтр хую чрьъэєЄюёЄш яЁхфяюыюцшь,
ўЄю яЁш $n \rightarrow \infty \ : h_n \rightarrow \ h $ \ ш \ $ W h_n \ \rightarrow \ g \ $.
%╥ръшь юсЁрчюь,
%$$
%\forall \  \varepsilon > 0, \ \exists \ n_0, \ \forall \ n \ \geq \ n_0 : \  \Vert \ W h_n \ - \ g \Vert < \ \varepsilon,
%$$
%Єю хёЄ№ фы  ы■сюую $  n \ \geq \ n_0 $  шьххь
╧юёыхфэхх ючэрўрхЄ, ўЄю \ $ f(x)h_n \ \rightarrow \ g(x) \ $ яю эюЁьх яЁюёЄЁрэёЄтр $ L^1 (X, \mathcal A, \mu, G)$ , Єю хёЄ№
$$
\lim\limits_{n\to\infty} \int\limits_X  \Vert \ f(x)h_n \ - \ g(x) \Vert \ d\mu(x) = 0.
$$
═ю Єюуфр ёє∙хёЄтєхЄ яюфяюёыхфютрЄхы№эюёЄ№ $ \{ h_{n_k} \}_{k=1}^{\infty}$ яюёыхфютрЄхы№эюёЄш $ \{ h_n \}_{n=1}^{\infty}$
Єрър , ўЄю \ $ f(x)h_{n_k} \ \rightarrow \ g(x) \ \ \mu $-яюўЄш тё■фє эр $ X $. ╬Єё■фр ёыхфєхЄ, ўЄю \ $ f(x)h = g(x) \ \mu $-яюўЄш тё■фє эр $X$, Єю
хёЄ№ \ $ W h \ = \ g \ $. ╥ръшь юсЁрчюь, юяхЁрЄюЁ $W$  ты хЄё  чрьъэєЄ√ь, р, чэрўшЄ, ш юуЁрэшўхээ√ь.

╟рьхЄшь, ўЄю $ T = SW $, уфх $ S :  L^1 (X, \mathcal A, \mu, G) \rightarrow G $ - юЄюсЁрцхэшх, яюЁюцфрхьюх шэЄхуЁрыюь ┴юїэхЁр. ╥ръ ъръ $ S $  ты хЄё  ёцрЄшхь, Єю
$ T $- юуЁрэшўхээ√щ ышэхщэ√щ юяхЁрЄюЁ,
% уфх $ g \ = \ g(x) \ \in \ L^1 (X, \mathcal A, \mu, G)$. ╥юуфр шч ыхьь√ ╘рЄє шьххь
%$$
%\Vert \ W h \ - \ g \Vert \ = \ \int\limits_X  \Vert \ f(x)h \ - \ g(x) \Vert \ d\mu(x) \ =
%$$
%$$
%= \ \int\limits_X  \ \ \underline{\lim}  \ \ \Vert \ f(x)h_n \ - \ g(x) \Vert \ d\mu(x) \  \leq \ \underline{\lim} \ \int\limits_X  \Vert \ f(x)h_n \ - \ g(x) \Vert \ %d\mu(x) \ < \ \varepsilon.
%$$
%┬ ёшыє яЁюшчтюы№эюёЄш $\varepsilon \ > \ 0,$ яюыєўрхь \ $ W h \ = \ g \ $, ш юяхЁрЄюЁ $W$  ты хЄё  чрьъэєЄ√ь, р, чэрўшЄ, ш юуЁрэшўхээ√ь.
%╥ръшь юсЁрчюь, фы  ы■сюую \ $h \ \in \ H $
%$$
%\int\limits_X \ \Vert \ f(x)h \Vert \ d\mu(x) \ =  \  \Vert \ W h \Vert  \leq \  \Vert \ W \ \Vert \ \Vert h \Vert .
%$$
%╬Єё■фр ёыхфєхЄ
%$$
%\Vert \ T h \Vert \ = \ \Vert \int\limits_X \ f(x)h \ d\mu(x) \Vert \ \leq \int\limits_X \ \Vert \ f(x)h \Vert \ d\mu(x) \ = \ \Vert \ W \ \Vert \ \Vert h \Vert,
%$$
ш ЄхюЁхьр фюърчрэр.

╬яЁхфхы хь√щ яюёыхфэхщ ЄхюЁхьющ юяхЁрЄюЁ эрч√трхЄё  ёшы№э√ь
шэЄхуЁрыюь ┴юїэхЁр, ъюЄюЁ√щ юсючэрўр■Є $(s) \int\limits_X fd\mu$.

╧єёЄ№ $( X, \mathcal A)$ - шчьхЁшьюх яЁюёЄЁрэёЄтю, $E$ - срэрїютю
яЁюёЄЁрэёЄтю. ╘єэъЎш  $\nu:\mathcal A\to E$ эрч√трхЄё  $E$-чэрўэ√ь чрЁ фюь эр
$(X, \mathcal A)$, хёыш фы  ы■сюую фшч·■эъЄэюую \ ёхьхщёЄтр $\{
A_n\}_{n=1}^\infty \subset A$
\begin{equation}
\nu \left ( \bigcup_{n=1}^\infty A_n \right ) = \sum_{n=1}^\infty
\nu (A_n).\label{21}
\end{equation}

┬ ёыєўрх $E={\mathbb C}$ ($E={\mathbb R}$) уютюЁ Є ю ъюьяыхъёэюь
(ъюэхўэюь тх∙хёЄтхээюь) чрЁ фх. ┬ ърўхёЄтх $E$ ьюцэю сЁрЄ№ ш
Ёрё°шЁхээюх ьэюцхёЄтю тх∙хёЄтхээ√ї ўшёхы $\bar{\mathbb R}$, ЁрчЁх°р 
 $\nu$ яЁшэшьрЄ№ ёрьюх сюы№°хх юфэю шч фтєї эхёюсёЄтхээ√ї чэрўхэшщ
$-\infty$ шыш $+\infty$.

┼ёыш $E={\c B} (H,G)$ уфх $H$ ш $G$ -- ушы№схЁЄют√ яЁюёЄЁрэёЄтр,
Єю чрЁ ф $\nu$ сєфхь эрч√трЄ№ ЁртэюьхЁэ√ь, хёыш ёїюфшьюёЄ№
Ё фр (\ref{21}) ЁрёёьрЄЁштрхЄё  яю юяхЁрЄюЁэющ эюЁьх, ш ёшы№э√ь,
хёыш ёїюфшьюёЄ№ юяхЁрЄюЁэюую Ё фр т (\ref{21}) ёшы№эр .

┬рЁшрЎшхщ чрЁ фр $\nu:\, {\mathcal A}\to E$ эрч√трхЄё  ЇєэъЎш 
$\vert \nu\vert :\, {\mathcal A}\to [0,+\infty]$, юяЁхфхы хьр 
ЇюЁьєыющ
$$
\vert \nu\vert (A) := {\rm sup} \sum_{n=1}^m \Vert \nu(A_n)\Vert,
\quad A\in {\mathcal A},
$$
уфх ёєяЁхьєь схЁхЄё  яю тёхь ъюэхўэ√ь Ёрчсшхэш ь $A$ эр фшч·■эъЄэ√х
яюфьэюцхёЄтр $\{A_n\}_{n=1}^m \subset {\mathcal A}$. ┬рЁшрЎш 
$\vert\nu\vert$ чрЁ фр $\nu$  ты хЄё  эршьхэ№°хщ яюыюцшЄхы№эющ ьхЁющ
 эр ${\mathcal A}$, ьрцюЁшЁє■∙хщ ЇєэъЎш■ $\Vert \nu(A) \Vert$, $A\in
 {\mathcal A}$.

┬ё ъшщ ъюэхўэ√щ тх∙хёЄтхээ√щ чрЁ ф (р, чэрўшЄ, ш ъюьяыхъёэ√щ) шьххЄ
ъюэхўэє■ трЁшрЎш■. ╬Єё■фр ёыхфєхЄ, ўЄю яЁш ${\rm dim} E<\infty$
трЁшрЎш  $\vert\nu\vert$ чрЁ фр $\nu$ ъюэхўэр. ╬фэръю, т ёыєўрх
схёъюэхўэюьхЁэюую срэрїютр яЁюёЄЁрэёЄтр $E$ ьюцхЄ юърчрЄ№ё , ўЄю
$\vert\nu\vert(X) = +\infty$. ┼ёыш $E={\c B} (H,G)$ уфх $H$
ш $G$ -- срэрїют√ яЁюёЄЁрэёЄтр, Єю ёЇюЁьєышЁютрээюх т√°х юяЁхфхыхэшх
трЁшрЎшш яЁшьхэшью ъръ ъ ЁртэюьхЁэ√ь, Єръ ш ъ ёшы№э√ь ${\c B}
(H,G)$-чэрўэ√ь чрЁ фрь. ╧Ёш ¤Єюь ёшы№э√щ чрЁ ф ё ъюэхўэющ
трЁшрЎшхщ  ты хЄё  ЁртэюьхЁэ√ь.

$E$-чэрўэ√щ чрЁ ф $\nu:\, {\mathcal A}\to E$ эрч√трхЄё  рсёюы■Єэю
эхяЁхЁ√тэ√ь юЄэюёшЄхы№эю ьхЁ√ $\mu:\, {\mathcal A} \to [0,+\infty]$,
хёыш фы  ы■сюую $A\in {\mathcal A}$ шч $\mu(A)=0$ т√ЄхърхЄ $\nu(A)=
0$. ┬рцэхщ°шь яЁшьхЁюь рсёюы■эю эхяЁхЁ√тэюую $E$-чэрўэюую чрЁ фр
(юЄэюёшЄхы№эю ьхЁ√ $\mu$)  ты хЄё  шэЄхуЁры ┴юїэхЁр
\begin{equation}
\nu(A): =\int\limits_A f(x) d\mu(x), \quad A\in {\mathcal A},
\label{22}
\end{equation}
уфх $f\in {\mathcal L}^1 (X, {\mathcal A},\mu, E)$. ╠эюцхёЄтю
чэрўхэшщ Єръюую чрЁ фр ёхярЁрсхы№эю. ╧Ёш ¤Єюь т ёыєўрх
$E=\bar{\mathbb R}$ шыш $E={\mathbb C}$ чрЁ ф $\nu$ шьххЄ ъюэхўэє■
трЁшрЎш■, ъюЄюЁр  т√ЁрцрхЄё  шэЄхуЁрыюь ╦хсхур
\begin{equation}
\vert\nu\vert (A) =\int\limits_A \vert f(x)\vert d\mu(x),\quad A\
\in {\mathcal A}. \label{23}
\end{equation}
─ы  ёъры Ёэ√ї чрЁ фют шьххЄ ьхёЄю ёыхфє■∙р 

\begin{theorem}
$($╤ъры Ёэ√щ трЁшрэЄ ЄхюЁхь√ ╨рфюэр-═шъюфшьр, $\cite{1}, \cite{6} \ )$
╧єёЄ№ $\nu:\, {\mathcal A} \to {\mathbb R} ({\mathbb C})$ --
чрЁ ф, рсёюы■Єэю эхяЁхЁ√тэ√щ юЄэюёшЄхы№эю $\sigma$-ъюэхўэющ ьхЁ√
$\mu:\, {\mathcal A} \to [0,+\infty ]$. ╥юуфр ёє∙хёЄтєхЄ
хфшэёЄтхээр  $($ ё ЄюўэюёЄ№■ фю ЁртхэёЄтр $\mu$-яюўЄш тё■фє $)$ ЇєэъЎш 
$f\in {\mathcal L}^1 (X,{\mathcal A},\mu) \ (f\in {\mathcal L}^1
(X,{\mathcal A},\mu, {\mathbb C})$) Єрър , ўЄю
$$
\nu(A)=\int\limits_A f(x) d\mu(x), \quad A\in {\mathcal A}.
$$
\end{theorem}

╘єэъЎш  $f$ шч ╥хюЁхь√ 2.6 эрч√трхЄё  яЁюшчтюфэющ ╨рфюэр-═шъюфшьр
чрЁ фр $\nu$ яю ьхЁх $\mu$, шыш яыюЄэюёЄ№■ чрЁ фр $\nu$ юЄэюёшЄхы№эю
ьхЁ√ $\mu$, ш юсючэрўрхЄё  $\frac{d\nu}{d\mu}$.

┬ юс∙хь ёыєўрх фы  $E$-чэрўэ√ї чрЁ фют т ёЇюЁьєышЁютрээюь тшфх
єЄтхЁцфхэшх ╥хюЁхь√ 2.6 эхтхЁэю. ╥ръ, эряЁшьхЁ, юэю эх шьххЄ ьхёЄр
фрцх
фы  ёыєўр  ъюэхўэющ ьхЁ√ $\mu$ ш ёхярЁрсхы№эюую ушы№схЁЄютюую
яЁюёЄЁрэёЄтр $E$. ╧Ёштхфхээ√щ эшцх яЁшьхЁ яюыєўхэ эхёє∙хёЄтхээющ
ьюфшЇшърЎшхщ
эр ёыєўрщ ушы№схЁЄютюую яЁюёЄЁрэёЄтр рэрыюушўэюую яЁшьхЁр т \cite{1}
 (Appendix E, Exercise 7).

%\begin{example}
╧ЁшьхЁ. ╧єёЄ№ $\lambda$--ьхЁр ╦хсхур эр $[0,1]$, ${\mathcal A}$ -- $\sigma$-
рыухсЁр шчьхЁшь√ї яю ╦хсхує яюфьэюцхёЄт юЄЁхчър $X=[0,1]$,
$E=L^2[0,1]$ ш чрЁ ф $\nu:\, {\mathcal A}\to E$ чрфрэ ЇюЁьєыющ
$\nu(A)=\chi_A$, $A\in {\mathcal A}$, уфх $\chi_A$ --
їрЁръЄхЁшёЄшўхёър  ЇєэъЎш  ьэюцхёЄтр $A$. ╤ўхЄэр  рффшЄштэюёЄ№
чрЁ фр $\nu$ яЁютхЁ хЄё  эхяюёЁхфёЄтхээю. └сёюы■Єэр  эхяЁхЁ√тэюёЄ№
$\nu$
юЄэюёшЄхы№эю $\lambda$ т√ЄхърхЄ шч ЁртхэёЄтр
$$
\Vert\nu(A)\Vert = \left ( \int\limits_0^1 \chi_A^2 (x) d\lambda (x)
\right )^{1/2} = \sqrt{\lambda(A)},\quad A\in {\mathcal A}.
$$
╬Єё■фр ыхуъю ёыхфєхЄ, ўЄю $\vert\nu\vert(X)=+\infty$. ╤ фЁєующ
ёЄюЁюэ√, хёыш с√ ёє∙хёЄтютрыр Єрър  ЇєэъЎш  $f\in
{\mathcal L}^1 (X,{\mathcal A},\mu, E)$, ўЄю т√яюыэ ыюё№ с√
(\ref{22}), Єю шч юяЁхфхыхэш  трЁшрЎшш чрЁ фр ёыхфютрыр с√ юЎхэър
$$
\vert\nu\vert (X) \le \int\limits_X \Vert f(x)\Vert dx <+\infty,
$$
ўЄю яЁштхыю с√ ъ яЁюЄштюЁхўш■.
%\end{example}

╩ръ ь√ єтшфшь фрыхх, фюяюыэшЄхы№эюх ЄЁхсютрэшх ъюэхўэюёЄш трЁшрЎшш
чрЁ фр фрхЄ тючьюцэюёЄ№ ЁрёяЁюёЄЁрэшЄ№ ЄхюЁхьє ╨рфюэр-═шъюфшьр
эр ёыєўрщ тхъЄюЁэючэрўэ√ї чрЁ фют.

\section{┬хъЄюЁэ√щ ш юяхЁрЄюЁэ√щ трЁшрэЄ√ ЄхюЁхь√ ╨рфюэр-═шъюфшьр}
\setcounter{equation}{0}

 ╤эрўрыр фюърцхь рэрыюуш яЁхфёЄртыхэш 
(\ref{23}) фы  трЁшрЎшш чрЁ фют, чрфрээ√ї шэЄхуЁрырьш ┴юїэхЁр.

\begin{lemma}
$(\cite{10})$ ╧єёЄ№ $(X,{\mathcal A}, \mu)$ -- яЁюёЄЁрэёЄтю ё $\sigma$
-ъюэхўэющ ьхЁющ, $E$ - ушы№схЁЄютю яЁюёЄЁрэёЄтю ш $f\in
{\mathcal L}^1 (X,{\mathcal A},\mu, E)$. ╥юуфр чрЁ ф
$$
\nu(A) : = \int\limits_A f(x) d\mu(x), \quad A\in {\mathcal A},
$$
шьххЄ ъюэхўэє■ трЁшрЎш■, ш
\begin{equation}
\vert\nu\vert(A) =\int\limits_A \Vert f(x)\Vert d\mu(x), \quad A\in
 {\mathcal A}.
\label{31}
\end{equation}
\end{lemma}

{\bf ─юърчрЄхы№ёЄтю.} ═х эрЁє°р  юс∙эюёЄш, ьюцэю ёўшЄрЄ№, ўЄю $E$
ёхярЁрсхы№эю. ┬ яЁюЄштэюь ёыєўрх, чрьхэшт ЇєэъЎш■ $f$ эр ёютярфр■∙є■
ё эхщ $\mu$-яюўЄш тё■фє ЇєэъЎш■ $f_1\in {\mathcal L}^1 (X,{\mathcal A}
,\mu, E)$ ё ёхярЁрсхы№э√ь ьэюцхёЄтюь чэрўхэшщ, ьюцэю ёєчшЄ№ $E$ фю
чрь√ърэш  ышэхщэющ юсюыюўъш ¤Єюую ьэюцхёЄтр.

╚ч эхЁртхэёЄтр
$$
\Vert\nu (A)\Vert \le \int\limits_A \Vert f(x)\Vert d\mu(x),
\quad A\in {\mathcal A},
$$
т√ЄхърхЄ юЎхэър
\begin{equation}
\vert\nu\vert (A) \le \int\limits_A \Vert f(x)\Vert d\mu(x), \quad
A\in {\mathcal A}. \label{32}
\end{equation}
╬Єё■фр ёыхфєхЄ ъюэхўэюёЄ№ трЁшрЎшш $\vert\nu\vert$ ш хх рсёюы■Єэр 
эхяЁхЁ√тэюёЄ№ юЄэюёшЄхы№эю ьхЁ√ $\mu$. ╤юуырёэю ╥хюЁхьх 2.6 ёє∙хёЄтєхЄ
эхюЄЁшЎрЄхы№эр  ЇєэъЎш  $p(x)\in {\mathcal L}^1 (X,{\mathcal A},
\mu)$ Єрър , ўЄю
\begin{equation}
\vert\nu\vert (A) = \int\limits_A p(x)d\mu(x), \quad A\in {\mathcal A}.
\label{33}
\end{equation}
╚ч (\ref{32}) ш (\ref{33}) ёыхфєхЄ, ўЄю фы  $\mu$-яюўЄш тёхї $x\in X$
\begin{equation}
p(x)\le\Vert f(x)\Vert.
\label{34}
\end{equation}
╤ фЁєующ ёЄюЁюэ√, яЁш ы■сюь $e\in E$ фы  ёъры Ёэюую чрЁ фр
$$
\eta_e(A):=<\nu(A),e>=\int\limits_A <f(x),e>d\mu(x), \quad A\in
{\mathcal A},
$$
яюыєўрхь
\begin{equation}
\vert\eta_e(A)\vert \le \Vert \nu(A)\Vert\Vert e\Vert \le
\vert\nu\vert(A)\Vert e\Vert =\int\limits_A (p(x) \Vert e\Vert)
d\mu(x), \quad A\in {\mathcal A}. \label{35}
\end{equation}
╤ыхфютрЄхы№эю, яЁш тёхї $e\in E$ шч (\ref{23}) ш {\ref{35}) т√ЄхърхЄ
\begin{equation}
\int\limits_A \vert <f(x), e>\vert d\mu(x) = \vert\eta_e\vert(A)
\le \vert\nu\vert(A)\Vert e\Vert =\int\limits_A (p(x)\Vert e\Vert)
d\mu(x), \quad A\in {\mathcal A}. \label{36}
\end{equation}
┬√сЁрт ёўхЄэюх тё■фє яыюЄэюх т $E$ яюфьэюцхёЄтю $E_1$, шч (\ref{36})
 яюыєўрхь, ўЄю ёє∙хёЄтєхЄ $A_0\in {\mathcal A}$ Єръюх, ўЄю
$\mu(A_0) =0$ ш яЁш ы■сюь $x\in X\setminus A_0$ ш ы■сюь $e\in E_1$
т√яюыэ хЄё  эхЁртхэёЄтю
$$
\vert <f(x),e>\vert \le p(x)\Vert e\Vert.
$$
╧хЁхїюф  т ¤Єюь эхЁртхэёЄтх ъ яЁхфхыє яЁш $e\to f(x)$, $e\in E_1$,
яюыєўрхь, ўЄю фы  $x\in X\setminus A_0$ (Є.х. $\mu$-яюўЄш
тё■фє эр $X$)
\begin{equation}
\Vert f(x)\Vert \le p(x).
\label{37}
\end{equation}
╚ч (\ref{33}), (\ref{34}) ш (\ref{37}) т√ЄхърхЄ (\ref{31}).

\begin{lemma}
╧єёЄ№ $(X,{\mathcal A}, \mu)$ -- яЁюёЄЁрэёЄтю ё $\sigma$-ъюэхўэющ
ьхЁющ, $H$ ш $G$ -- ушы№схЁЄют√ яЁюёЄЁрэёЄтр, яЁшўхь $H$ ёхярЁрсхы№эю ш
$f:\, X\to {\c B} (H,G)$ -- ёшы№эю шэЄхуЁшЁєхьр  яю ┴юїэхЁє
ЇєэъЎш . ─ы  Єюую ўЄюс√ ёшы№э√щ чрЁ ф
$$
\nu (A) : = (s)\int\limits_A f(x) d\mu(x),\quad A\in {\mathcal A},
$$
шьхы ъюэхўэє■ трЁшрЎш■, эхюсїюфшью ш фюёЄрЄюўэю, ўЄюс√ $\Vert f(x)
\Vert \in {\mathcal L}^1 (X, {\mathcal A},\mu)$. ╧Ёш ¤Єюь
\begin{equation}
\vert\nu\vert (A) =  \int\limits_A \Vert f(x)\Vert d\mu(x),
\quad A\in {\mathcal A}. \label{38}
\end{equation}
\end{lemma}

{\bf ─юърчрЄхы№ёЄтю.} ╧Ёхцфх тёхую яюърцхь, ўЄю шч $f(x) h\in
{\mathcal L}^1 (X, {\mathcal A}, \mu, G)$ яЁш тёхї $h\in H$ т√ЄхърхЄ
шчьхЁшьюёЄ№ ЇєэъЎшш $\Vert f(x) \Vert$. ─хщёЄтшЄхы№эю, т√сЁрт
ёўхЄэюх яыюЄэюх эр хфшэшўэющ ёЇхЁх т $H$ яюфьэюцхёЄтю $S = \{ h_n\vert
n\in {\mathbb N}\}$, яюыєўрхь яЁш тёхї $x\in X$
\begin{equation}
\Vert f(x)\Vert = \sup\limits_n \Vert f(x) h_n\Vert,
\label{39}
\end{equation}
юЄъєфр ёыхфєхЄ шчьхЁшьюёЄ№ $\Vert f(x) \Vert$.

╧єёЄ№ трЁшрЎш  чрЁ фр $\nu$ ъюэхўэр. ╥юуфр т ёшыє хую рсёюы■Єэющ
эхяЁхЁ√тэюёЄш юЄэюёшЄхы№эю ьхЁ√ $\mu$ яю ╥хюЁхьх 2.6 яюыєўрхь
яЁхфёЄртыхэшх тшфр (\ref{33}), уфх $p(x)\in {\mathcal L}^1 (X,
{\mathcal A}, \mu)$ -- эхюЄЁшЎрЄхы№эр  ЇєэъЎш . ╥юуфр яЁш ы■сюь
$h\in H$ фы  $G$-чэрўэюую чрЁ фр
$$
\lambda_h(A):=\nu(A)h=\int\limits_A (f(x)h)d\mu(x),\quad A\in
{\mathcal A},
$$
т√яюыэ хЄё 
\begin{equation}
\Vert\lambda_h (A)\Vert \le \Vert \nu(A)\Vert\Vert h\Vert\le
\vert\nu\vert (A)\Vert h\Vert =\int\limits_A (p(x)\Vert h\Vert)d\mu(x),
\quad A\in {\mathcal A}.
\label{310}
\end{equation}
╤ыхфютрЄхы№эю, трЁшрЎш  чрЁ фр $\lambda_h$ ъюэхўэр. ╚ч (\ref{310})
т ёшыє ╦хьь√ 3.1 т√ЄхърхЄ
\begin{equation}
\int\limits_A \Vert f(x)h\Vert d\mu(x) =\vert \lambda_h\vert
(A)\le \vert\nu\vert(A)\Vert h\Vert =\int\limits_A (p(x)\Vert
h\Vert)d\mu(x), \quad A\in {\mathcal A}. \label{311}
\end{equation}
╚ч (\ref{311}) ёыхфєхЄ, ўЄю ёє∙хёЄтєхЄ $A_0 \in {\mathcal A}$ Єръюх,
 ўЄю $\mu(A_0)=0$ ш яЁш ы■сюь $x\in X\setminus A_0$ ш яЁш ы■сюь
$h_n\in S$ т√яюыэ хЄё  юЎхэър
\begin{equation}
\Vert f(x) h_n\Vert \le p(x) \Vert h_n\Vert = p(x),
\label{312}
\end{equation}
Єръ ъръ $\Vert h_n\Vert =1$. ╚ч (\ref{39}) ш (\ref{312}) т√ЄхърхЄ,
ўЄю яЁш $\mu$-яюўЄш тёхї $x\in X$
\begin{equation}
\Vert f(x)\Vert \le p(x),
\label{313}
\end{equation}
ш, ёыхфютрЄхы№эю, $\Vert f(x)\Vert \in {\mathcal L}^1 (X,
{\mathcal A}, \mu)$.

╬сЁрЄэю, яєёЄ№ $\Vert f(x)\Vert \in {\mathcal L}^1 (X,{\mathcal A},
\mu)$. ╥юуфр шч ёшы№эющ шэЄхуЁшЁєхьюёЄш яю ┴юїэхЁє ЇєэъЎшш $f$
яЁш ы■сюь $h\in H$ яюыєўрхь
$$
\Vert \nu(A)h\Vert \le \int\limits_A \Vert f(x)h\Vert d\mu(x) \le
\int\limits_A (\Vert f(x)\Vert \Vert h\Vert) d\mu (x),
\quad A\in {\mathcal A}.
$$
╬Єё■фр т√ЄхърхЄ юЎхэър
$$
\Vert \nu(A)\Vert \le \int\limits_A \Vert f(x)\Vert d\mu(x), \quad A\in
{\mathcal A}.
$$
╤ыхфютрЄхы№эю, трЁшрЎш  чрЁ фр $\nu$ ъюэхўэр. ╚ч
ъюэхўэюёЄш трЁшрЎшш $\vert\nu\vert$ ш хх рсёюы■Єэющ
эхяЁхЁ√тэюёЄш ёыхфєхЄ яЁхфёЄртыхэшх (\ref{33}).
╥хяхЁ№, ъръ ш Ёрэхх, яюыєўрхь эхЁртхэёЄтр (\ref{34}) ш (\ref{313}), шч
ъюЄюЁ√ї т√ЄхърхЄ (\ref{38}).

\begin{rem}
╬ЄьхЄшь, ўЄю яЁш т√яюыэхэшш єёыютшщ ╦хьь√ 3.2 фрцх т ёыєўрх
ъюэхўэюёЄш трЁшрЎшш чрЁ фр $\nu$ ёшы№э√щ шэЄхуЁры ┴юїэхЁр эх
 ты хЄё , тююс∙х уютюЁ , ЁртэюьхЁэ√ь, їюЄ  ёрь чрЁ ф $\nu$ ЁртэюьхЁхэ.
\end{rem}

─хщёЄтшЄхы№эю, яєёЄ№, эряЁшьхЁ, $X=[0,1]$, ${\mathcal A}$ --
$\sigma$-рыухсЁр шчьхЁшь√ї яю ╦хсхує яюфьэюцхёЄт юЄЁхчър $[0,1]$,
$\lambda$ -- ьхЁр ╦хсхур эр $[0,1]$, $H=G=L^2(X,{\mathcal A},
\lambda)$. ╨рёёьюЄЁшь $f:\, X\to{\c B} (H)$, уфх $f(x)$,
$x\in [0,1],$ -- юЁЄюяЁюхъЄюЁ т $H$, чрфртрхь√щ єьэюцхэшхь эр
їрЁръЄхЁшёЄшўхёъє■ ЇєэъЎш■ $\chi_{[0,x]}$ юЄЁхчър $[0,x]$. ╘єэъЎш 
$f$, юўхтшфэю, ёшы№эю шэЄхуЁшЁєхьр яю ┴юїэхЁє, ш т ёшыє $\Vert
f(x)\Vert =1$ ($0<x\le 1$) ёююЄтхЄёЄтє■∙шщ чрЁ ф $\nu$ шьххЄ
ъюэхўэє■ трЁшрЎш■. ┬ьхёЄх ё Єхь ЇєэъЎш  $f$ эх  ты хЄё 
$\mu$-ёє∙хёЄтхээю ёхярЁрсхы№эющ
 т ${\c B}(H)$, Єръ ъръ яЁш $x_i\in X$ ($i=1,2$) ш
$x_1\not= x_2$ шьххь $\Vert f(x_1) - f(x_2)\Vert=1$. ┴юыхх Єюую,
$f$ эх  ты хЄё  $({\mathcal A}, {\mathcal B})$-шчьхЁшьющ, уфх
${\mathcal B}$ -- $\sigma$-рыухсЁр сюЁхыхтёъшї ьэюцхёЄт т срэрїютюь
яЁюёЄЁрэёЄтх ${\c B}(H)$. ┬ ёрьюь фхых, шч фшёъЁхЄэюёЄш
Єюяюыюушш, шэфєЎшЁютрээющ эр $f(X)$, ш шэ·хъЄштэюёЄш ЇєэъЎшш $f$
т√ЄхърхЄ, ўЄю тёх яюфьэюцхёЄтр ьэюцхёЄтр $[0,1]$  ты ■Єё 
яЁююсЁрчрьш эхъюЄюЁ√ї юЄъЁ√Є√ї ьэюцхёЄт шч ${\mathcal B}$.

\begin{theorem}
$($┬хъЄюЁэ√щ трЁшрэЄ ЄхюЁхь√ ╨рфюэр-═шъюфшьр.$)$ ╧єёЄ№ $(X,{\mathcal A},
\mu)$ -- яЁюёЄЁрэёЄтю ё $\sigma$-ъюэхўэющ ьхЁющ, $E$ -- ушы№схЁЄютю
яЁюёЄЁрэёЄтю ш $\nu:\, {\mathcal A}\to E$ -- чрЁ ф, єфютыхЄтюЁ ■∙шщ
єёыютш ь:

$($a$)$ ьэюцхёЄтю чэрўхэшщ чрЁ фр $\nu$ ёхярЁрсхы№эю;

$($с$)$ чрЁ ф $\nu$ рсёюы■Єэю эхяЁхЁ√тхэ юЄэюёшЄхы№эю ьхЁ√ $\mu$;

$($т$)$ трЁшрЎш  чрЁ фр $\nu$ ъюэхўэр.

╥юуфр ёє∙хёЄтєхЄ хфшэёЄтхээр  $($ё ЄюўэюёЄ№■ фю ЁртхэёЄтр $\mu$-яюўЄш
тё■фє$)$ ЇєэъЎш  $f\in {\mathcal L}^1 (X, {\mathcal A}, \mu, E)$
Єрър , ўЄю
$$
\nu(A) = \int\limits_A f(x) d\mu(x), \quad A\in {\mathcal A}.
$$
\end{theorem}

{\bf ─юърчрЄхы№ёЄтю.} ╩ръ ш т ╦хььх 3.1, эх юуЁрэшўштр  юс∙эюёЄш,
яЁюёЄЁрэёЄтю $E$ ьюцэю ёўшЄрЄ№ ёхярЁрсхы№э√ь. ┬√схЁхь т $E$
ъръющ-эшсєф№
юЁЄюэюЁьшЁютрээ√щ срчшё $\{e_k \}_{k=1}^\infty$. ╥юуфр
\begin{equation}
\nu(A) =\sum_{n=1}^\infty <\nu(A), e_k> e_k, \quad A\in {\mathcal A}.
\label{314}
\end{equation}
╤ъры Ёэ√х чрЁ ф√ $\eta_k(A):=<\nu(A), e_k>$, $A\in {\mathcal A}$,
$k\in{\mathbb N}$, єфютыхЄтюЁ ■Є эхЁртхэёЄтє
$\vert \eta_k(A)\vert \le \Vert \nu(A)\Vert$. ╤ыхфютрЄхы№эю, фы 
шї трЁшрЎшщ $\vert\eta_k\vert$  т√яюыэ ■Єё  юЎхэъш
$$
\vert\eta_k\vert(A)\le \vert\nu\vert(A), \quad A\in {\mathcal A}, \ \
k\in{\mathbb N}.
$$
╥ръшь юсЁрчюь, чрЁ ф√ $\eta_k$, $k\in{\mathbb N}$, рсёюы■Єэю
эхяЁхЁ√тэ√ юЄэюёшЄхы№эю $\mu$. ╧ю ╥хюЁхьх 2.6 яЁш $k\in{\mathbb N}$
ёє∙хёЄтє■Є $q_k \in {\mathcal L}^1 (X, {\mathcal A}, \mu)$ Єрър ,
ўЄю
$$
\eta_k (A) = \int\limits_A q_k (x) d\mu(x), \quad A\in {\mathcal A}.
$$
╬ўхтшфэю, ўЄю яЁш ы■сюь $m\in{\mathbb N}$ ЇєэъЎш 
\begin{equation}
f_m (x) := \sum_{k=1}^m q_k (x) e_k
\label{315}
\end{equation}
яЁшэрфыхцшЄ ${\mathcal L}^1 (X, {\mathcal A}, \mu, E)$. ╬ЄьхЄшь,
ўЄю $E$-чэрўэ√х чрЁ ф√
$$
\nu_m(A) = \sum_{k=1}^m <\nu(A),e_k>e_k, \ \ A\in {\mathcal A}, \  m\in{\mathbb N},
$$
т ёшыє ышэхщэюёЄш шэЄхуЁрыр ┴юїэхЁр фюяєёър■Є яЁхфёЄртыхэшх
\begin{equation}
\nu_m (A) = \int\limits_A f_m (x) d\mu(x), \quad A\in {\mathcal A}, \
m\in{\mathbb N}.
\label{316}
\end{equation}
╚ч (\ref{314}) тшфэю, ўЄю яЁш ы■с√ї $A\in {\mathcal A}$ ш
$m\in{\mathbb N}$ чрЁ ф $\nu_m$ єфютыхЄтюЁ хЄ эхЁртхэёЄтє $\Vert
\nu_m(A)\Vert \le \Vert \nu(A)\Vert$, \ р, чэрўшЄ, ш эхЁртхэёЄтє
\begin{equation}
\vert\nu_m\vert(A)\le \vert \nu\vert(A), \quad A\in {\mathcal A}, \
m\in{\mathbb N}.
\label{317}
\end{equation}
╚ч (\ref{316}) ёюуырёэю ╦хььх 3.1 яюыєўрхь, ўЄю
\begin{equation}
\vert\nu_m\vert (A)=\int\limits_A \Vert f_m (x)\Vert d\mu(x),
\quad A\in {\mathcal A}, \  m\in{\mathbb N}.
\label{318}
\end{equation}
╩Ёюьх Єюую, яю ╥хюЁхьх 2.6 ёє∙хёЄтєхЄ эхюЄЁшЎрЄхы№эр  ЇєэъЎш 
$p\in {\mathcal L}^1 (X, {\mathcal A}, \mu)$ Єрър , ўЄю
\begin{equation}
\vert \nu\vert (A) =\int\limits_A p(x) d\mu(x), \quad A\in
{\mathcal A}.
\label{319}
\end{equation}
╚ч (\ref{317}) -- (\ref{319}) т√ЄхърхЄ, ўЄю яЁш ы■сюь
$m\in{\mathbb N}$ фы  $\mu$-яюўЄш тёхї $x\in X$ т√яюыэ хЄё 
эхЁртхэёЄтю
\begin{equation}
\Vert f_m (x)\Vert \le p(x).
\label{320}
\end{equation}
╬Єё■фр ш шч (\ref{315}) ёыхфєхЄ, ўЄю фы  $\mu$-яюўЄш тёхї
$x\in X$
$$
\sum_{k=1}^\infty \vert q_k (x)\vert^2 \le p^2 (x).
$$

╤ыхфютрЄхы№эю, фы  $\mu$-яюўЄш тёхї $x\in X$ ёїюфшЄё  Ё ф
$\sum_{k=1}^\infty q_k (x) e_k$, Єю хёЄ№ яЁш $\mu$-яюўЄш тёхї $x\in X$
юяЁхфхыхэр ЇєэъЎш  $f(x):= \lim\limits_{m\to\infty} f_m (x)$.
╧юырур  $f(x) : = 0$ эр шёъы■ўшЄхы№эюь ьэюцхёЄтх, яюыєўрхь, ўЄю т ёшыє ╥хюЁхь√ 2.2
ЇєэъЎш   $f$ ёшы№эю $\mu$-шэЄхуЁшЁєхьр.
╬Єё■фр, єўшЄ√тр  (\ref{316}) ш (\ref{320}), яю ╥хюЁхьх 2.4
яюыєўрхь, ўЄю
$f\in {\mathcal L}^1 (X, {\mathcal A}, \mu, E)$ ш
$$
\nu(A) = \int\limits_A f(x) d\mu(x), \quad A\in {\mathcal A}.
$$
┼ёыш $\hat f(x)\in {\mathcal L}^1 (X, {\mathcal A}, \mu, E)$ --
фЁєур  ЇєэъЎш  Єрър , ўЄю
$$
\nu(A) = \int\limits_A \hat f(x) d\mu(x), \quad A\in {\mathcal A},
$$
Єю $\int\limits_A (f(x) -\hat f(x)) d\mu(x)=0$,
$A\in {\mathcal A}$. ═ръюэхЎ,  ёюуырёэю  ╦хььх 3.1  шьххь  
$$ 
\int\limits_A \Vert f(x) -\hat f(x)\Vert d\mu(x)=0, \ \ \  A\in {\mathcal A},
$$
 ш, чэрўшЄ, $f(x) = \hat f(x)$
 яЁш $\mu$-яюўЄш тёхї $x\in X$.

\begin{theorem}
$($╬яхЁрЄюЁэ√щ трЁшрэЄ ЄхюЁхь√ ╨рфюэр-═шъюфшьр.$)$ ╧єёЄ№ $(X,
{\mathcal A}, \mu)$ -- яЁюёЄЁрэёЄтю ё $\sigma$-ъюэхўэющ ьхЁющ, $H$
ш $G$ -- ушы№схЁЄют√
яЁюёЄЁрэёЄтр, яЁшўхь $H$ ёхярЁрсхы№эю, ш $\nu:\, {\mathcal A}\to
{\c B}(H,G)$ -- ЁртэюьхЁэ√щ чрЁ ф, єфютыхЄтюЁ ■∙шщ єёыютш ь:

$($р$)$  яЁш ы■сюь $h\in H$ \ \ $G$-чэрўэ√щ чрЁ ф $\lambda_h (A):= \nu(A) h$,
 $A\in {\mathcal A}$, шьххЄ ёхярЁрсхы№эюх ьэюцхёЄтю чэрўхэшщ;

$($с$)$ чрЁ ф $\nu$ рсёюы■Єэю эхяЁхЁ√тхэ юЄэюёшЄхы№эю ьхЁ√ $\mu$;

$($т$)$ трЁшрЎш  чрЁ фр $\nu$ ъюэхўэр.

╥юуфр ёє∙хёЄтєхЄ хфшэёЄтхээр  $($ё ЄюўэюёЄ№■ фю ЁртхэёЄтр $\mu$-яюўЄш
тё■фє$)$ ёшы№эю шэЄхуЁшЁєхьр  яю ┴юїэхЁє ЇєэъЎш  $f:\, X\to
{\c B}(H,G)$
Єрър , ўЄю
\begin{equation}
 \nu(A) = (s)\int\limits_A f(x) d\mu(x), \quad A\in {\mathcal A}.
\label{321}
\end{equation}
\end{theorem}

{\bf ─юърчрЄхы№ёЄтю.} ╧Ёш ы■сюь $h\in H$ чрЁ ф $\lambda_h$
єфютыхЄтюЁ хЄ эхЁртхэёЄтє $\Vert \lambda_h(A)\Vert \le \Vert
\nu(A)\Vert \Vert h\Vert$,
$A\in {\mathcal A}$. ╬Єё■фр фы  трЁшрЎшш шьххь юЎхэъє
\begin{equation}
\vert\lambda_h\vert(A) \le \vert\nu\vert(A) \Vert h\Vert,\quad A\in
{\mathcal A}, \ \ h\in H.
\label{322}
\end{equation}
╥ръшь юсЁрчюь, тёх чрЁ ф√ $\lambda_h$, $h\in H$ рсёюы■Єэю
эхяЁхЁ√тэ√ юЄэюёшЄхы№эю $\mu$ ш шьх■Є ъюэхўэє■ трЁшрЎш■. ╧ю
╥хюЁхьх 3.4
яЁш ы■сюь $h\in H$ ёє∙хёЄтєхЄ хфшэёЄтхээр  (ё ЄюўэюёЄ№■ фю
ЁртхэёЄтр $\mu$-яюўЄш тё■фє) ЇєэъЎш  $f_h\in {\mathcal L}^1
(X, {\mathcal A}, \mu, G)$
Єрър , ўЄю
\begin{equation}
\nu(A)h =\lambda_h (A) = \int\limits_A f_h (x) d\mu(x), \quad A\in
{\mathcal A}.
\label{323}
\end{equation}
╧юърцхь, ўЄю ЇєэъЎш■ $f_h (x)$, $h\in H$, ьюцэю т√сЁрЄ№ Єръ, ўЄюс√
яЁш ърцфюь $x\in X$ чртшёшьюёЄ№ юЄ $h$ с√ыр ышэхщэющ, Єю хёЄ№,
ўЄюс√ т√яюыэ ыюё№ ЁртхэёЄтю
\begin{equation}
f_h (x) = f(x) h, \quad  h\in H, x\in X,
\label{324}
\end{equation}
уфх $f:\, X\to {\c B} (H,G)$. ╥юуфр шч (\ref{323}),
(\ref{324}), ╥хюЁхь√ 2.5 ш юяЁхфхыхэш  $\lambda_h$ сєфхЄ т√ЄхърЄ№
(\ref{321}).

╠√ сєфхь ёўшЄрЄ№ ушы№схЁЄютю яЁюёЄЁрэёЄтю $H$ схёъюэхўэюьхЁэ√ь, эх
юёЄрэртыштр ё№ эр яЁюёЄюь ёыєўрх яЁюёЄЁрэёЄтр ъюэхўэющ ЁрчьхЁэюёЄш.
┬√схЁхь яЁюшчтюы№э√щ юЁЄюэюЁьшЁютрээ√щ срчшё
$\{h_k\}_{k=1}^\infty$ т $H$. ╨рёёьюЄЁшь ышэхщэє■ юсюыюўъє $H_1$
тхъЄюЁют $\{h_k\}_{k=1}^\infty$ эрф яюыхь ${\mathbb Q}$ т ёыєўрх
тх∙хёЄтхээ√ї $H$ ш $G$ шыш эрф яюыхь ${\mathbb Q}+i{\mathbb Q}$ т
ёыєўрх ъюьяыхъёэ√ї $H$ ш $G$. ▀ёэю, ўЄю $H_1$ -- ёўхЄэюх тё■фє
яыюЄэюх яюфьэюцхёЄтю т $H$. ╟рфрфшь яЁш ърцфюь $x\in X$ юяхЁрЄюЁ $f(x)$
%чэрўхэшх ЇєэъЎшш $f(x)$, $x\in X$,
ёэрўрыр эр тхъЄюЁрї срчшёр $\{h_k\}_{k=1}^\infty$, яюыюцшт
$f(x)h_k := f_{h_k} (x)$, $k\in{\mathbb N}$. ╥юуфр фы  ы■сюую
$k\in{\mathbb N}$ шьххь $f(x)h_k \in {\mathcal L}^1 (X, {\mathcal
A}, \mu,G)$ ш т ёшыє (\ref{323}) яюыєўрхь ЁртхэёЄтю
$$
\nu(A)h_k = \int\limits_A (f(x)h_k) d\mu(x), \quad A\in
{\mathcal A}.
$$

─рыхх яЁш ърцфюь $x\in X$ юяЁхфхышь  юяхЁрЄюЁ $f(x)$ эр ы■сюь $h=\sum_{k=1}^m
\alpha_k h_k \in H_1$, яюыюцшт $f(x) h:= \sum_{k=1}^m \alpha_k f(x)
h_k$. ╥хь
ёрь√ь яЁш тёхї $h\in H_1$ ЇєэъЎш  $f(x) h\in {\mathcal L}^1 (X,
{\mathcal A}, \mu,G),$ ш яЁш ¤Єюь  фы  ы■сюую $A\in {\mathcal A}$ ёяЁртхфышт√ ЁртхэёЄтр
\begin{equation}
\nu(A)h = \sum_{k=1}^m \alpha_k \nu(A)h_k = \sum_{k=1}^m \alpha_k \int\limits_A (f(x)h_k))  = \int\limits_A (  \sum_{k=1}^m \alpha_k f(x)
h_k )d\mu(x) = \int\limits_A (f(x)h) d\mu(x).
\label{325}
\end{equation}
╚ч єёыютшщ ({\it с}), ({\it т}) ш ╥хюЁхь√ 2.6 яюыєўрхь, ўЄю ёє∙хёЄтєхЄ
эхюЄЁшЎрЄхы№эр  ЇєэъЎш  $p(x)\in {\mathcal L}^1 (X, {\mathcal A},
\mu)$ Єрър , ўЄю
\begin{equation}
\vert \nu\vert (A) = \int\limits_A p(x)d\mu(x), \quad A\in
{\mathcal A}. \label{326}
\end{equation}
╚ч (\ref{322}), (\ref{323}), (\ref{325}), (\ref{326}) ш ╦хьь√ 3.1 т√ЄхърхЄ,
ўЄю яЁш тёхї $h\in H_1$
$$
\int\limits_A \Vert f(x)h\Vert d\mu(x)\le \int\limits_A (p(x)
\Vert h\Vert) d\mu(x),  \quad A\in {\mathcal A},
$$
юЄъєфр яЁш ърцфюь $h\in H_1$ фы  $\mu$-яюўЄш тёхї $x\in X$
ёяЁртхфыштю эхЁртхэёЄтю
\begin{equation}
\Vert f(x)h\Vert \le p(x)\Vert h\Vert.
\label{327}
\end{equation}
┬ ёшыє ёўхЄэюёЄш $H_1$ ёє∙хёЄтєхЄ Єръюх $A_0\in A$, ўЄю $\mu
(A_0)=0$ ш яЁш ¤Єюь яЁш ы■с√ї $x\in X\setminus A_0$ \ ш \ $h\in H_1$
шьххЄ ьхёЄю юЎхэър (\ref{327}), Єю хёЄ№ ышэхщэ√щ юяхЁрЄюЁ $f(x)$ юуЁрэшўхэ
эр $H_1$ (ъръ юяхЁрЄюЁ фхщёЄтє■∙шщ шч $H_1$ т $G$ эрф яюыхь
${\mathbb Q}$ шыш ${\mathbb Q}+i{\mathbb Q}$ ёююЄтхЄёЄтхээю). ▌Єю
яючтюы хЄ яЁш $x\in X\setminus A_0$ яЁюфюыцшЄ№ яю эхяЁхЁ√тэюёЄш
 юяхЁрЄюЁ $f(x)$ эр $H$ (ъръ юяхЁрЄюЁ фхщёЄтє■∙шщ шч $H$ т $G$ эрф
яюыхь ${\mathbb R}$ шыш ${\mathbb C}$ ёююЄтхЄёЄтхээю). └ шьхээю,
хёыш $h\in H$ ш $\{k_n\}_{n=1}^\infty\subset H_1$ Єрър , ўЄю
$h=\lim\limits_{n\to\infty} k_n$, Єю яЁш тёхї $x\in X\setminus A_0$
 яюыюцшь
\begin{equation}
f(x)h:=\lim\limits_{n\to\infty} (f(x) k_n).
\label{328}
\end{equation}
╬ўхтшфэю, Єръюх яЁюфюыцхэшх ъюЁЁхъЄэю. ╥ръшь юсЁрчюь, яЁш
$x\in X\setminus A_0$ юяхЁрЄюЁ $f(x) \in {\c B}(H,G)$ ш
єфютыхЄтюЁ хЄ єёыютш■
\begin{equation}
\Vert f(x)h\Vert \le p(x) \Vert h\Vert, \quad h\in H.
\label{329}
\end{equation}
╧Ёш $x\in A_0$ яхЁхюяЁхфхышь $f(x)$ эр тхъЄюЁрї $h\in H_1$,
яюыюцшт $f(x)h=0$, ш яЁюфюыцшь $f(x)$ эєыхь эр тёх $H$. ╥хь ёрь√ь
яЁш ы■сюь $x\in X$
юяхЁрЄюЁ $f(x)\in {\c B}(H,G)$, яЁш ¤Єюь ЁртхэёЄтю (\ref{328})
 ш юЎхэър (\ref{329}) т√яюыэ ■Єё  яЁш тёхї $x\in X$. ╬Єё■фр, єўшЄ√тр , ўЄю  $f(x) k_n
\in {\mathcal L}^1 (X, {\mathcal A}, \mu,G)$ ш яюёыхфютрЄхы№эюёЄ№ $ \Vert k_n \Vert $ юуЁрэшўхэр,
ёюуырёэю ╥хюЁхьрь 2.2 ш 2.4
яюыєўрхь, ўЄю яЁш ы■сюь $h\in H$ ЇєэъЎш  $f(x)h \in {\mathcal L}^1
(X, {\mathcal A}, \mu,G)$
ш ёяЁртхфыштю ЁртхэёЄтю
$$
\nu(A)h = \int\limits_A (f(x)h) d\mu(x), \quad A\in
{\mathcal A},
$$
Єю хёЄ№ т√яюыэ хЄё  ЁртхэёЄтю (\ref{323}) яЁш $f_h(x) :=f(x)h$.

─юърчрЄхы№ёЄтю хфшэёЄтхээюёЄш $f$ яЁютюфшЄё  Єръ цх, ъръ ш т
╥хюЁхьх 3.4, яЁш ¤Єюь тьхёЄю ╦хьь√ 3.1 эхюсїюфшью тюёяюы№чютрЄ№ё 
╦хььющ 3.2.

\begin{cor}
┼ёыш т ╥хюЁхьх 3.5 чрЁ ф $\nu$ тьхёЄю єёыютш  $(${\it р}$)$
єфютыхЄтюЁ хЄ сюыхх ёшы№эюьє ЄЁхсютрэш■:

$(\dot{a}) \  \nu$ шьххЄ ёхярЁрсхы№эюх ьэюцхёЄтю чэрўхэшщ,

Єю ЇєэъЎш  $f(x)$  ты хЄё  ЁртэюьхЁэю шэЄхуЁшЁєхьющ яю ┴юїэхЁє ш
шэЄхуЁры $($\ref{321}$)$  ты хЄё  ЁртэюьхЁэ√ь.
\end{cor}

{\bf ─юърчрЄхы№ёЄтю. } ┬ ¤Єюь ёыєўрх тьхёЄю срэрїютр яЁюёЄЁрэёЄтр
${\c B}(H,G)$ ьюцэю ЁрёёьюЄЁхЄ№ хую эршьхэ№°хх чрьъэєЄюх
яюфяЁюёЄЁрэёЄтю $F$,
ёюфхЁцр∙хх ьэюцхёЄтю чэрўхэшщ чрЁ фр $\nu$. ╬эю, юўхтшфэю,  ты хЄё 
ёхярЁрсхы№э√ь срэрїют√ь яЁюёЄЁрэёЄтюь.

╧ЁютхЁшь $({\mathcal A}, {\mathcal B})$-шчьхЁшьюёЄ№ юЄюсЁрцхэш  $f$,
 уфх ${\mathcal B}$ -- $\sigma$-рыухсЁр сюЁхыхтёъшї ьэюцхёЄт т $F$.
 ─ы 
¤Єюую т ёшыє ёхярЁрсхы№эюёЄш $F$ фюёЄрЄюўэю яЁютхЁшЄ№, ўЄю яЁш ы■сюь
$r> 0$ ш ы■сюь $T\in F$ т√яюыэ хЄё  $\{ x\in X: \Vert f(x) -
T\Vert \le r\}
\in {\mathcal A}$. ╧юёыхфэхх т√ЄхърхЄ шч ЁртхэёЄтр
$$
\{ x\in X: \Vert f(x) - T\Vert \le r\}= \bigcap\limits_{h\in S}
\{ x\in X: \Vert f(x)h - Th\Vert \le r\},
$$
уфх $S$ -- ёўхЄэюх тё■фє яыюЄэюх ьэюцхёЄтю эр хфшэшўэющ ёЇхЁх т $H$,
 ш ёшы№эющ шчьхЁшьюёЄш $G$-чэрўэ√ї ЇєэъЎшщ $f(x)h$ яЁш ърцфюь $h\in H$.

{\underline{╩юььхэЄрЁшщ}}. ┬ ЁрсюЄх \cite{14} ўрёЄэ√щ ёыєўрщ
╥хюЁхь√ 3.4 ЇюЁьєышЁєхЄё  фы  сюЁхыхтёъюую чрЁ фр эр хфшэшўэющ
юъЁєцэюёЄш ${\mathbb T}$,
эю яЁш ¤Єюь (╥хюЁхьр 3.2, ёЄЁ. 111) юЄёєЄёЄтєхЄ єёыютшх (т) \ ╥хюЁхь√ 3.4 фрээющ ЁрсюЄ√.
%ртЄюЁ√ ю°шсюўэю ёўшЄрыш, ўЄю
%єёыютшх (b)  ты хЄё  ёыхфёЄтшхь єёыютш  (с).
└эрыюушўэр  эхЄюўэюёЄ№
 шьхыр ьхёЄю ш
ё юяхЁрЄюЁэючэрўэ√ьш ьхЁрьш (╦хььр 3.3, ёЄЁ.114), уфх, ъЁюьх
 Єюую, юЄёєЄёЄтютрыю єёыютшх ({\it р}) ╥хюЁхь√ 3.5 фрээющ ЁрсюЄ√.
 ╬фэръю чрьхЄшь, ўЄю тёх Ёхчєы№ЄрЄ√ ЁрсюЄ√ \cite{14}, яЁш фюърчрЄхы№ёЄтх
 ъюЄюЁ√ї шёяюы№чютрышё№ ёююЄтхЄёЄтхээю ╥хюЁхьр 3.2 (ёЄЁ. 111) ш ╦хььр 3.3 (ёЄЁ. 114)
 юёЄр■Єё  т ёшых, Єръ ъръ т ¤Єшї яЁшыюцхэш ї т√°хєяюь эєЄюх єёыютшх (т), р т юяхЁрЄюЁэюь ёыєўрх ш
 єёыютшх (р), т√яюы ■Єё  ртЄюьрЄшўхёъш.

╬сэрЁєцхээ√х ртЄюЁрьш (єцх яюёых яєсышърЎшш) эхЄюўэюёЄш яюсєфшыш шї
яЁюрэрыш-чшЁютрЄ№ ышЄхЁрЄєЁє яю фрээ√ь тюяЁюёрь. ╧Ёш ¤Єюь, ъръ єцх
юЄьхўхэю
тю ┬тхфхэшш, эх Єюы№ъю эх єфрыюё№ эрщЄш яЁшхьыхьюх шчыюцхэшх
ЁрёёьюЄЁхээ√ї т фрээющ ёЄрЄ№х трЁшрэЄют ЄхюЁхь√ ╨рфюэр-═шъюфшьр,
эю ш с√ыш
юсэрЁєцхэ√ х∙х эхъюЄюЁ√х яюуЁх°эюёЄш т ЇюЁьєышЁютърї, яхЁхфр■∙шхё 
юЄ юфэшї ртЄюЁют ъ фЁєушь. ╥ръ, эряЁшьхЁ, т ьюэюуЁрЇшш \cite{13}
╥хюЁхьр 1.2 (ёЄЁ.325) ёЇюЁьєышЁютрэр т сюыхх юс∙хщ ЇюЁьх, ўхь
╥хюЁхьр 3.5 фрээющ ЁрсюЄ√. ┬ ¤Єющ ЄхюЁхьх (яЁш шёяюы№чєхь√ї эрьш
юсючэрўхэш ї)
$X={\mathbb R}$, ${\mathcal A}$ -- $\sigma$-рыухсЁр сюЁхыхтёъшї
ьэюцхёЄт эр ${\mathbb R}$, $H$ -- ёхярЁрсхы№эюх срэрїютю
яЁюёЄЁрэёЄтю, $G$ --
ёюяЁ цхээюх яЁюёЄЁрэёЄтю рэЄшышэхщэ√ї эхяЁхЁ√тэ√ї ЇєэъЎшюэрыют эр
эхъюЄюЁюь ёхярЁрсхы№эюь срэрїютюь яЁюёЄЁрэёЄтх, $\mu$ -- трЁшрЎш 
ёырсю
$\sigma$-рффшЄштэюую юяхЁрЄюЁэючэрўэюую чрЁ фр $\nu:\,
{\mathcal A}\to {\c B} (H,G)$, ъюЄюЁр  ыюъры№эю юуЁрэшўхэр.
╧Ёш ¤Єюь юЄёєЄёЄтєхЄ
єёыютшх $(\dot{a})$ (ёь. ╤ыхфёЄтшх 3.6 фрээющ ЁрсюЄ√), эю єЄтхЁцфрхЄё , ўЄю ёырс√щ
шэЄхуЁры, ъюЄюЁ√ь яЁхфёЄртыхэ чрЁ ф $\nu$ эр ьхЁх $\mu$  ты хЄё 
Єръцх
ёїюф ∙шьё  яю юяхЁрЄюЁэющ эюЁьх эр ьэюцхёЄтрї ъюэхўэющ ьхЁ√ (Є.х.
яЁхфёЄртшь эр ¤Єшї ьэюцхёЄтрї ЁртэюьхЁэ√ь шэЄхуЁрыюь ┴юїэхЁр).
╬ЄёєЄёЄтшх т ╥хюЁхьрї 1.1 ш 1.2 т \cite{13} фюърчрЄхы№ёЄтр ёшы№эющ
шчьхЁшьюёЄш яыюЄэюёЄш $ \Psi: {\mathbb R} \to {\c B} (H,G) $ чрЁ фр $\nu$
эх яючтюы хЄ уютюЁшЄ№ ю ёє∙хёЄтютрэшш шэЄхуЁрыр $\int\limits_{\Delta} \Psi(\lambda) d\mu(\lambda)$
эр ьэюцхёЄтрї ъюэхўэющ ьхЁ√. ╬фэръю, єЄтхЁцфхэшх ╥хюЁхь√ 1.1 т \cite{13} тхЁэю т ёшыє ёхярЁрсхы№эюёЄш
ьэюцхёЄтр чэрўхэшщ юяхЁрЄюЁэющ ьхЁ√. ▌Єю ёыхфєхЄ, эряЁшьхЁ, шч ╤ыхфёЄтш  3.6 фрээющ ЁрсюЄ√. ╫Єю цх ърёрхЄё 
єЄтхЁцфхэш  ╥хюЁхь√ 1.2 т \cite{13}, Єю
ъюэЄЁяЁшьхЁюь ъ эхьє ьюцхЄ ёыєцшЄ№ чрЁ ф $\nu$, ЁрёёьюЄЁхээ√щ т
╟рьхўрэшш 3.3 ъ ╦хььх 3.2 фрээющ ЁрсюЄ√, хёыш тьхёЄю $X = [0,1]$ тч Є№ $X = {\mathbb R}$ , р $f(x), \ x\in {\mathbb R},$
юяЁхфхышЄ№ ъръ юяхЁрЄюЁ т $H$, чрфртрхь√щ єьэюцхэшхь эр їрЁръЄхЁшёЄшўхёъє■ ЇєэъЎш■ $\chi_{[-\infty,x]}$ ыєўр $[-\infty,x].$
▌ЄюЄ чрЁ ф яЁхфёЄртшь ёшы№э√ь шэЄхуЁрыюь ┴юїэхЁр (р,чэрўшЄ, ш ёырс√ь шэЄхуЁрыюь) эр ьэюцхёЄтрї ъюэхўэющ ьхЁ√, эю
эх яЁхфёЄртшь ЁртэюьхЁэ√ь шэЄхуЁрыюь ┴юїэхЁр.

%яюыєўхээ√щ юўхтшфэющ
%ьюфшЇшърЎшхщ эр ёыєўрщ $X={\mathbb R}$ чрЁ фр, ЁрёёьюЄЁхээюую т
%╟рьхўрэшш 3.3
%ъ ╦хььх 3.2 фрээющ ЁрсюЄ√. ▌ЄюЄ чрЁ ф яЁхфёЄртшь ёшы№э√ь шэЄхуЁрыюь ┴юїэхЁр (р,
%чэрўшЄ, ш ёырс√ь шэЄхуЁрыюь), эю эх яЁхфёЄртшь ЁртэюьхЁэ√ь
%шэЄхуЁрыюь ┴юїэхЁр.
%═хтхЁэю, ёфхырээюх т \cite{13} яЁхфяюыюцхэшх, ўЄю фюърчрЄхы№ёЄтю
%╥хюЁхь√ 1.2 "ьюцэю яЁютюфшЄ№ яю ёїхьх фюърчрЄхы№ёЄтр ╥хюЁхь√ 1.1"
%(ёЄЁ. 325),
%яюёъюы№ъє єЄтхЁцфхэшх яюёыхфэхщ ЄхюЁхь√ тхЁэю т ёшыє
%ёхярЁрсхы№эюёЄш ьэюцхёЄтр чэрўхэшщ юяхЁрЄюЁэющ ьхЁ√.
┬ фры№эхщ°хь ¤Єр ЄхюЁхьр, єцх ъръ ЄхюЁхьр
┴хЁхчрэёъюую-├хы№Їрэфр-╩юёЄ■ўхэъю, ЎшЄшЁєхЄё  ш ЇюЁьєышЁєхЄё 
ё ЁртэюьхЁэ√ь шэЄхуЁрыюь
┴юїэхЁр (ёь.,
эряЁшьхЁ, \cite{15}, ╧Ёхфыюцхэшх 2.15, ёЄЁ.19). ▌Єю юсёЄю Єхы№ёЄтю  тшыюё№
фюяюыэшЄхы№э√ь ьюЄштюь эряшёрэш  фрээющ ЁрсюЄ√.

└тЄюЁ√ сыруюфрЁ Є ┬.╠.╩рфхЎр, тэшьрЄхы№эю яЁюўшЄрт°хую
яхЁтюэрўры№э√щ трЁшрэЄ ЁрсюЄ√ ш т√ёърчрт°хую Ё ф яюыхчэ√ї чрьхўрэшщ.

\end{document}